\documentclass[11pt, reqno, a4paper]{amsart}
\usepackage[margin=1.2in]{geometry}
\numberwithin{equation}{section}
\usepackage{amssymb,amsfonts,amsthm}
\usepackage[utf8]{inputenc}
\usepackage{listings}
\usepackage{bm}
\usepackage[normalem]{ulem}
\usepackage{xcolor}
\definecolor{myred}{RGB}{200, 30, 30}

\definecolor{myblue}{RGB}{200, 30, 30}

\usepackage{color}
\usepackage[bookmarks]{hyperref}

\allowdisplaybreaks[4]

\newtheoremstyle{myremark}{10pt}{10pt}{}{}{\bfseries}{.}{.5em}{}

 \newtheorem{thm}{Theorem}[section]
 \newtheorem{cor}[thm]{Corollary}
 \newtheorem{lem}[thm]{Lemma}
 \newtheorem{prop}[thm]{Proposition}
 \theoremstyle{definition}
 \newtheorem{defn}[thm]{Definition}
 
 \newtheorem{rem}[thm]{Remark}

\newcommand{\R}{\ensuremath{\mathbb R}}

\allowdisplaybreaks[4]

\begin{document}

\title[Stability of the Morrey and Hardy inequalities]{Gradient Stability of the Morrey and Hardy inequalities}

\keywords{Morrey inequality, Hardy inequality,  quantitative gradient stability, Schwarz symmetrization, capacity.}

\subjclass{46E35, 26D10, 49J40.}

\author[U. Das]{ Ujjal Das}
\address[Ujjal Das]{School of Mathematical Sciences, National Institute of Science Education and Research, Bhubaneswar 
752050, India.}
\email{udas@niser.ac.in, getujjaldas@gmail.com}

\author[S. Roy]{ Subhajit Roy}
\address[Subhajit Roy]{Department of Mathematics, Indian Institute of Science Education and Research, Pune 411008, India.}
\email{ rsubhajit.math@gmail.com, subhajit.roy@iiserpune.ac.in}

\thanks{}



\begin{abstract} 
We find that the quantitative gradient stability estimate for the Morrey inequality in a ball $B_R$ fails with respect to any exponent of the {\it gradient-level} distance $D_M(u):=\inf\{\|\nabla u-\nabla v\|_{L^p(B_R)}:v\in \mathcal{M}\}$. Instead, the estimate holds with $D_M(u^\star),$ the gradient-level distance between the Schwarz symmetrization $u^\star$ and the space of Morrey extremals $\mathcal{M}$. Nevertheless, we establish a gradient stability estimate for the Morrey inequality by slightly `weakening' the above distance. We also obtain analogous gradient stability for the geometric $L^p$-Hardy inequality in certain domains. Distinct behaviors in the stability of these two inequalities have been discussed. Furthermore, in the $L^2$-case, we investigate the optimal constants and their attainability in both the Morrey and Hardy gradient stability inequalities.

\end{abstract}
\maketitle
\date{March 2026}

\section{Introduction}
In recent years, the study of {\it quantitative stability} of various functional inequalities has received significant attention. 
Let
$\mathcal{E}(u)\geq \mathcal{F}(u)$
be a functional inequality on a function space $X$, where $\mathcal{E}(\cdot), \mathcal{F}(\cdot)$ are certain energy functionals. 
Assume that the set of extremal functions 
$\mathcal{M}:=\{u\in X:\mathcal{E}(u)=\mathcal{F}(u)\}$ is non-empty. A fundamental question in the theory of functional inequalities is whether or not the inequality admits a \emph{quantitative stability estimate}. More precisely, the question is: if $u$ is \textit{almost an extremal}, 
  i.e., the {\it deficit functional}  $\delta(u):=\mathcal{E}(u)-\mathcal{F}(u)$ is close to $0$, then is it necessary that 
 \(u\) has to be `close' to the set of extremals \(\mathcal{M}\)? Such a property is typically expressed by a quantitative estimate of the form
\begin{align} \label{Eq:gen_stab}
    \mathcal{E}(u)-\mathcal{F}(u)\geq C\,d(u,\mathcal{M}) \ \ \ \forall u \in X,
\end{align}
where $C>0$ is a constant and $d(u,\mathcal{M})$ is a suitable notion of distance from $u$ to $\mathcal{M}$.
The question of quantitative stability for the $L^2$-Sobolev inequality was first raised by Brezis and Lieb \cite{Brezis-Lieb} and was settled by Bianchi and Egnell \cite{Bianchi}. 
They showed that
there exists $C>0$ such that for every $u\in W^{1,2}_0(\mathbb{R}^N),$
\begin{align*} 
     \|\nabla u\|_{L^2(\mathbb{R}^N)}^2-\mathcal{S}^2\|u\|_{L^{2^*}(\mathbb{R}^N)}^2\geq C\,\inf_{v\in\mathcal{M}}\|\nabla u-\nabla v\|_{L^2(\mathbb{R}^N)}^2,
\end{align*} where {$N \geq 3$,} $2^*=\frac{2N}{N-2}$, and $\mathcal{S}$ is the optimal Sobolev constant. The Hilbert space structure of $W^{1,2}_0(\mathbb{R}^N)$ was crucial in Bianchi--Egnell's approach. Using symmetrization techniques and {optimal transport}, the quantitative stability for the $L^p$-Sobolev inequality was established by Cianchi, Fusco, Maggi, and Pratelli \cite{Ci-Fus2009} with the following distance function on $W^{1,p}_0(\mathbb{R}^N)$:  $$\inf_{v\in\mathcal{M}}\left(\frac{\| u- v\|_{L^{p^*}(\mathbb{R}^N)}}{\|u\|_{L^{p^*}(\mathbb{R}^N)}}\right)^\alpha\,,\qquad  \;\; p \in (1, N),$$ for some $\alpha>0.$ Later, for $p\ge 2,$ Figalli and Neumayer \cite{ Fig-New2019} obtained a stronger {\it gradient stability estimate}, in the spirit of Bianchi--Egnell's result, with respect to the distance  given by \begin{align} \label{fig}
    \inf_{v\in\mathcal{M}}\left(\frac{\| \nabla u- \nabla v\|_{L^p(\mathbb{R}^N)}}{\|\nabla u\|_{L^{p}(\mathbb{R}^N)}}\right)^\beta\quad \mbox{on} \ \  W^{1,p}_0(\mathbb{R}^N) \,,
\end{align} where $\beta=p\alpha$ with $\alpha$ as above. Stability in terms of the above {\it gradient-level distance} is particularly significant, as it directly quantifies the deviation at the level of the Sobolev energy. Neumayer \cite{Neu2020} further extended the previous gradient stability to the full range $1 < p < N$, and more recently, Figalli and Zhang \cite{Fi-Zh2022} obtained the optimal exponent $\beta$ in \eqref{fig}. Quantitative stability of the Sobolev inequality was also established in the fractional setting \cite{Chen-F-W, Chen}, Hyperbolic spaces \cite{Bhakta}, and Riemannian manifolds \cite{Nobili}. {Moreover, such quantitative stability for various other functional inequalities has been obtained recently.}
The literature in this direction is vast, and the references below are by no means exhaustive. We mention only a few works on quantitative stability that might be relevant to the objectives of the present paper; e.g., on  Gagliardo-Nirenberg-Sobolev inequality \cite{carlen,Zhang-GNS}, Caffarelli-Kohn-Nirenberg inequality \cite{wei-caf, Anh-Lu}, log-Sobolev inequality \cite{Jean,Chen}, and Poincar\'e inequality \cite{Lu-poin,Amir}.  
{ Quantitative stability of the classical Hardy inequality with respect to the distance from its {\it virtual extremals} was established in \cite[Euclidean case]{Chianchi_2008} and \cite[Cartan-Hadamard manifold]{Prasun}.} 

A quantitative estimate like \eqref{Eq:gen_stab} naturally leads us to the question of finding/estimating the optimal constant in \eqref{Eq:gen_stab}, and subsequently, whether the optimal constant is attained or not. In the linear case $p=2,$ the attainment and lower bound of the {\it{best stability constant}} is studied for the Sobolev inequality \cite{ konig, Jean}, Hardy-Sobolev inequality \cite{Souptik}, and Caffarelli-Kohn-Nirenberg inequality \cite{Cazacu, wei-Wu}; see the references therein for further development on this topic. 

In most of the above results, the dimension of $\mathcal{M}$ is greater than $1$. In this article, we investigate the {quantitative} gradient stability of the following two fundamental inequalities in bounded domains for which  $\mathcal{M}$ is one-dimensional:
    \begin{itemize}
        \item {\it{Morrey inequality}} (the counterpart of the Sobolev inequality when the Sobolev exponent $p$ is strictly greater than the dimension $N$),
        \item {\it{geometric $L^p$-Hardy inequality}}.
    \end{itemize}

Let $\Omega\subset\mathbb{R}^N$ be a bounded domain and $1<p< \infty$. {It is well known that $\|u\|_{W^{1,p}_0(\Omega)}:=\|\nabla u\|_{L^p(\Omega)}$ defines a norm on $W^{1,p}_0(\Omega)$ due to the Poincar\'e inequality.} 
Clearly, $W^{1,2}_0(\Omega)$ is an inner product space with inner product $\langle u,v\rangle_{W^{1,2}_0(\Omega)}:=\int_\Omega \nabla u\cdot\nabla v\,\mathrm{d}x\,,$ for all $u,v\in W^{1,2}_0(\Omega).$

For $p>N\ge 1$, the {\it{Morrey inequality}} states that 
\begin{equation}\label{morrey}
   \|\nabla u\|_{L^p(\Omega)}\ge C_{\mathrm{M}}\|u\|_{L^\infty(\Omega)} \qquad\forall\, u\in W^{1,p}_0(\Omega)\,,
\end{equation}
where $C_{\mathrm{M}}>0$ is the {\it optimal} constant in the above inequality.

For $p\in (1,\infty),$ we say
that the \textit{geometric $L^p$-Hardy inequality} is satisfied in $\Omega$ if 
\begin{equation} \label{Hardy-p}
 \|\nabla u\|_{L^p(\Omega)}\ge C_{\mathrm{H}}\|u\|_{L^p(\Omega,\,d_\Omega^{-p})}\qquad\forall\, u\in W^{1,p}_0(\Omega),
\end{equation}
where $C_{\mathrm{H}}>0$ is the {\it optimal} constant and the distance of a point $x\in \Omega$ to the boundary $\partial \Omega$ is denoted by $$d_\Omega(x):=\mathrm{dist}(x,\partial \Omega) \quad \mbox{and} \  \ \ \|u\|_{L^p(\Omega,\,d_\Omega^{-p})}:=\left(\int_\Omega\frac{|u|^p}{d_\Omega^p}{\mathrm{d}}x\right)^{1/p}.$$  
It is known that if $\Omega$ is convex, then the {space/manifold of extremals}, $\mathcal{M}_{\mathrm{Morrey}}$, for the Morrey inequality \eqref{morrey} is one-dimensional
 \cite{Lindgren}. The same is true for the $L^p$-Hardy inequality \eqref{Hardy-p} for $C^{1,\gamma}$ bounded domains with {\it positive spectral gap} \cite[Theorem 4.1 \& 5.1]{Pinch19} (defined in Section \ref{sec-Hardy}).

In this article, we investigate the quantitative gradient stability of the Morrey and Hardy inequalities with respect to the corresponding distances $D_M(\cdot)$ and $D_H(\cdot)$, respectively, as defined in \eqref{dis} and \eqref{dis_H}. We encounter several different behaviors in the stability of these two inequalities. For the Morrey inequality, $\mathcal{M}_{\mathrm{Morrey}}$ is one-dimensional when $\Omega=(-R,R)\subset\mathbb{R}$, and in this case there do not exist constants $C,\gamma>0$ such that $
\delta_M(u)\geq C D_M(u)^\gamma
$ on $ W^{1,2}_0(\Omega)$ (Remark \ref{Rem1}). In contrast, the counterpart holds for the Hardy inequality (Theorem \ref{Hardy-Thm-2}); there exists $C>0$ such that
$\delta_H(u)\geq C D_H(u)^2
$ on $ W^{1,2}_0(\Omega)$, where $\Omega \subset \R^N$ is a bounded $C^{1,\gamma}$ domain with positive spectral gap, in which case $\mathcal{M}_{\mathrm{Hardy}}$ is one-dimensional.
 On the other hand, if $\Omega$ is a ball then the space of extremals for the Hardy inequality is trivial (as the spectral gap is $0$), and consequently, the gradient stability does not hold in this case, see Proposition \ref{Rem-Hardy}. For the Morrey inequality in a ball, our results  (Theorem \ref{Morr-th-1d} \& Theorem \ref{M:S:THM}) show that if the deficit $\delta_M(u)$ of the Morrey inequality is close to $0$, then instead of $u$, the {\it Schwarz symmetrization} $u^\star$ is close to $\mathcal{M}_{\mathrm{Morrey}}$ with respect to the distance $D_M(\cdot)$. We call this phenomenon as {\it gradient pseudo-stability} for the Morrey inequality in a ball. In fact, we obtain  quantitative gradient pseudo-stability inequalities \eqref{morrey-2} and \eqref{simplify} for the Morrey inequality on a ball.  We observe that the analogous  gradient pseudo-stability inequality does not hold for the Hardy inequality in a ball, see Proposition \ref{Rem-Hardy}.

Since there do not exist $C,\gamma>0$ such that $\delta_M(u)\ge CD_M(u)^\gamma$ on $W^{1,2}_0((-R,R))$, $R>0$, but the   pseudo-stability inequality $\delta_M(u) \ge D_M(u^\star)^2$ on $W^{1,2}_0((-R,R))$ holds (Theorem \ref{Morr-th-1d}), one may anticipate {\it weakening} the distance function $D_M(\cdot)$ `suitably' in order to establish the quantitative stability of the Morrey inequality in a ball. Indeed, it turns out that there exists $C>0$ such that
\[\delta_M(u) \geq C \min\{D_M(u)^{2},D_M(u)^{4}\} \ \ \ \quad \mbox{on} \ \ W^{1,2}_0((-R,R)) \,,\]
see Theorem \ref{M:S:THM} for a precise general higher-dimensional statement. On the other hand, we establish the corresponding quantitative gradient stability for the $L^p$-Hardy inequality in for $C^{1,\gamma}$ bounded domains with { positive spectral gap}, see Theorem \ref{Hardy-Thm-2} ($p=2$), Theorem \ref{Hardy-Thm-grad} ($p>2$), and Theorem \ref{Hardy-Thm-grad<2} ($p<2$). Furthermore, in the $L^2$-case, we investigate the {\it{optimal constants}} and their {\it{attainability}} in both the Morrey and Hardy gradient stability inequalities, see Theorem \ref{Morr-th-1d} and Theorem \ref{Hardy-Thm-2}. 

Our approaches to proving the gradient stability of the Morrey and Hardy inequalities are significantly different. In particular, an important {\it capacity estimate} (Proposition \ref{Prop:cap}) for singleton sets in a ball plays a crucial role in the proof of stability for the former inequality, whereas a {\it weighted Poincar\'e inequality} (Lemma \ref{poincare}) and the {\it simplified energy functional} (Lemma \ref{prop-1}) are essential ingredients in establishing stability for the latter. Nevertheless, standard reduction-type lemmas (Lemma \ref{Thm:Aux} and Lemma \ref{ThmHardy:Aux}) have been commonly used in both cases.

We organize the article as follows: 
in Section \ref{main-results}, we state our main results concerning the quantitative stability estimates for the Morrey and Hardy inequalities. Sections \ref{Morrey_section} and \ref{sec-Hardy} are devoted to the proofs of the stability results of Morrey (Theorems \ref{Morr-th-1d}, \ref{M:S:THM},  \ref{thm:MS:sup}) and Hardy (Theorems \ref{Hardy-Thm-2}, \ref{Hardy-Thm-grad}, \ref{Hardy-Thm-grad<2}) inequalities, respectively, as well as {we} derive several auxiliary estimates used in their proofs. 
In the Appendix, we establish gradient
pseudo-stability of the Morrey inequality in $\R^N$ (Theorem \ref{Rn-stability}) and some auxiliary results, including a regularity estimate.

\section{Main Results}\label{main-results}
\subsection{ Morrey inequality}\label{Morrey-result}
 In the first part of this article, we study the quantitative stability of the Morrey inequality \eqref{morrey}.
If $\Omega=\mathbb{R}^N$, Talenti \cite{Talenti1994}
showed  that the Morrey inequality for weakly differentiable function $u$ in $\R^N$ with $|\operatorname{supp}u|<\infty$ and $|\nabla u|\in L^p(\R^N)$:
\begin{equation}\label{Rn_morrey}
\kappa|\operatorname{supp}u|^{\frac1N-\frac1p}\|\nabla u\|_{L^p(\mathbb{R}^N)}\ge  \|u\|_{L^\infty(\mathbb{R}^N)}\,,  \qquad \ \  \kappa
:=
N^{-\frac{1}{p}}\omega_N^{-\frac{1}{N}}
\left(
\frac{p-1}{p-N}
\right)^{\frac{1}{p'}} 
\end{equation}
admits the following extremals
\[
v_{a,b,x_0}(x)
=
\begin{cases}
a\left(
b^{\frac{p-N}{p-1}}
-
|x-x_0|^{\frac{p-N}{p-1}}
\right)
& \text{if } |x-x_0|<b, \\[1ex]
0
& \text{if } |x-x_0|\ge b \,,
\end{cases} \qquad a \in \R,\,b>0, \  \mbox{and} \  x_0\in \mathbb{R}^N \,.
\]
Here, $\omega_N$ denotes the Lebesgue measure of the unit ball in $\mathbb{R}^N$.
 Cianchi \cite[Theorem 1.1]{cian2008} established  quantitative stability for the Morrey inequality in $\mathbb{R}^N$ \eqref{Rn_morrey} for $N \geq 2,$ {using machinery from symmetrization theory}. More precisely, there exist positive constants $\alpha$ and $C$ such that
\begin{equation}\label{Morrey_Rn}
 \frac{
\kappa|\operatorname{supp}u|^{\frac1N-\frac1p}\|\nabla u\|_{L^p(\mathbb R^N)}
-\|u\|_{L^\infty(\mathbb R^N)}}{
\|u\|_{L^\infty(\mathbb R^N)}
}
\ge C\left(\inf_{a,b,x_0}
\frac{\left\|
u-v_{a,b,x_0}
\right\|_{L^\infty(\mathbb R^N)}}{\|u\|_{L^\infty(\mathbb R^N)}}\right)^\alpha.   
\end{equation}
It remains open  whether or not the distance on the right-hand side of \eqref{Morrey_Rn} can be replaced by the following  stronger gradient-level distance:
\[ \inf_{a,b,x_0}
\frac{\left\|
\nabla u- \nabla v_{a,b,x_0}
\right\|_{L^p(\mathbb R^N)}}{\|  u\|_{L^\infty(\mathbb R^N)}} \,.\]

We start our investigation with bounded convex domains in $\R^N$. In \cite{Lindgren}, authors proved the  existence and simplicity of extremal functions for \eqref{morrey} if $\Omega$ is a convex domain. In particular, if $\Omega$ is a ball, then we know these extremals explicitly. Precisely, if $\Omega=B_R(x_0)$ then the space of extremal functions in the Morrey inequality \eqref{morrey} is given by  \cite{Talenti1994}
   $$\mathcal{M}_{\mathrm{Morrey}}=\left\{c
   \,\Phi: c\in \mathbb{R}\right\}, \quad \mbox{where} \ \ \  \Phi(x)=R^\frac{p-N}{p-1}-|x-x_0|^\frac{p-N}{p-1} \,. $$ 
      Moreover, the optimal constant in \eqref{morrey} is given by
 $$\mathbb{M}_{p}(\Omega)=\left(\frac{p-N}{p-1}\right)^\frac{1}{p'}R^\frac{N-p}{p}(N\omega_N)^{\frac{1}{p}}, \qquad  p'=\frac{p}{p-1} \,.$$ 
It turns out that the analogous quantitative stability as in \eqref{Morrey_Rn} cannot hold with respect to the { stronger gradient-level distance}  when the underlying domain is a ball. More precisely, the $p$-Morrey\textit{ deficit} functional
  \[
\delta_{M}(u):=
\begin{cases}
\frac{\| \nabla u\|_{L^p(\Omega)}^{p'}- \mathbb{M}_{p}(\Omega)^{p'}\|u\|^{p'}_{L^\infty(\Omega)}}{\|  u\|_{L^\infty(\Omega)}^{p'}}
& \text{if } p<2, \\[3ex]
\frac{\| \nabla u\|_{L^p(\Omega)}^{p}- \mathbb{M}_{p}(\Omega)^{p}\|u\|^{p}_{L^\infty(\Omega)}}{\| u\|_{L^\infty(\Omega)}^{p}} 
& \text{if } p\ge 2\,, 
\end{cases} \qquad \quad u\in W^{1,p}_0(\Omega)
\] 
can not be bounded from below by any positive exponent of the following distance defined by
\begin{align} \label{dis}
   D_M(u):=\inf\left\{\frac{\|\nabla u-\nabla v\|_{L^p(\Omega)}}{\| u\|_{L^\infty(\Omega)}}:v\in \mathcal{M}_{\mathrm{Morrey}}\right\} \,, \quad u\in W^{1,p}_0(\Omega) \,,
\end{align}
see Remark \ref{Rem1}. 
Nevertheless, we prove a sharp pseudo-stability inequality for the Morrey inequality in the one-dimensional linear case.
\begin{thm}\label{Morr-th-1d}
 Let   $\Omega$ be a symmetric interval centered at the origin and $p=2$. Then
    \begin{equation}\label{morrey-2}
       \delta_M(u) \ge D_M(u^\star)^2,\quad \forall\,u\in W^{1,2}_0(\Omega). 
    \end{equation}
    Moreover, inequality \eqref{morrey-2} is sharp, and equality is attained by non-trivial functions.
    \end{thm}
\begin{rem}\label{Rem1} Let $\Omega=(-R,R),\,p=2\,,$ and $a\in (0,R).$ Define
\[u_a(x):=
\begin{cases}
\frac{R+x}{R+a}
& \text{if } x\in (-R,a), \\[1ex]
\frac{R-x}{R-a}
& \text{if } x\in [a,R)\,. 
\end{cases} 
\] Then $\|u_a\|_{L^\infty(\Omega)}=1$ and $\|u_a'\|_{L^2(\Omega)}^2=\frac{2R}{R^2-a^2},$
which implies that $\delta_M(u_a)=\frac{2a^2}{R(R^2-a^2)}.$ 
It can be verified that $D_M(u_a)^2=4aR/[(R^2-a^2)(R+a)].$
Thus, for any $\gamma>0,$  the ratio $\frac{\delta_M(u_a)}{D_M(u_a)^\gamma} \to 0$ by taking $a\to 0$ or $a\to R$ depending on the range of $\gamma$. Hence, the gradient stability $\delta_M(u)\ge CD_M(u)^\gamma$ {{fails}}.
\end{rem}

The above one-dimensional result indicates that it is necessary to `weaken' the distance function $D_M(\cdot)$ suitably in order  to establish the gradient stability of the Morrey inequality. Towards this, we prove the following higher-dimensional quantitative stability of the Morrey inequality.
\begin{thm}\label{M:S:THM}
   Let $p>N\ge 1$ and $\Omega$ be a ball in $\mathbb{R}^N$. Then there exists a positive constant $C,$ depending
only on $p,\,N,$ and the radius of $\Omega$  such that for every $u\in W^{1,p}_0(\Omega),$
\begin{align} \label{simplify}
    \delta_M(u)  \ge C
\begin{cases} 
\min\{D_M(u)^{p'},D_M(u)^{2p'}\}
& \text{if } p<2\,,
 \\[2ex]
\min\{D_M(u)^{p},D_M(u)^{p^2}\}
& \text{if } p\ge 2\,.
\end{cases} 
\end{align}
\end{thm}
In particular, if $\Omega=B_R(0)$ then  we can replace $D_M(u)$ by $D_M(u^\star)$ on the right-hand side of \eqref{simplify} (due to P\'olya--Szeg\"o inequality \ref{polya}), which gives a quantitative gradient pseudo-stability inequality for the Morrey inequality.

Observe that the distance function $D_M(\cdot)$ given by \eqref{dis} is not scale invariant, i.e., $D_M(\cdot)$ changes under the transformation $\Omega \rightarrow t\Omega$ for $t\ne 1.$  For this reason, we consider the following scale-invariant gradient-level distance as in \eqref{fig}
\begin{equation*}
  \tilde{D}_M(u):=\inf\left\{\frac{\|\nabla u-\nabla v\|_{L^p(\Omega)}}{\| \nabla u\|_{L^p(\Omega)}}:v\in \mathcal{M}_{\mathrm{Morrey}}\right\} \,.  
\end{equation*}
It is not difficult to verify that $D_M(u)\ge\mathbb{M}_p(\Omega)\tilde{D}_M(u) \,.$ Hence,
as a corollary of the above theorem, we obtain the following analogue of the stability inequality \eqref{simplify}. 
\begin{cor}\label{Morrey-cor}
Let $p>N\ge 1$ and $\Omega$ be a ball in $\R^N$. Then there exists a positive constant $C,$ depending
only on $p,\,N,$ and the radius of $\Omega,$ such that for every $u\in W^{1,p}_0(\Omega),$
\[
\delta_M(u)\ge C
\begin{cases}
\tilde{D}_M(u)^{2p'}
& \text{if } p<2\,, \\[2ex]
\tilde{D}_M(u)^{p^2}
& \text{if } p\ge 2\,.
\end{cases}
\]
\end{cor}

To prove Theorem \ref{M:S:THM}, we first establish a stability inequality in terms of a {`weaker'} distance given by \eqref{eq:ci_dist} (see Theorem \ref{thm:MS:sup}). A careful analysis of the proof of Theorem \ref{Morr-th-1d} together with the example in Remark \ref{Rem1}, shows us the need for an important capacity estimate for singleton sets in a ball, which is derived in Proposition \ref{Prop:cap} and it plays a key role in proving Theorem \ref{M:S:THM}. We then establish a reduction lemma (Lemma  \ref{Thm:Aux}) that helps us to connect the `weakened' distance \eqref{eq:ci_dist} with the distance $D_M(\cdot)$. Combining the `weak' stability inequality in Theorem \ref{thm:MS:sup} and the reduction lemma in Lemma  \ref{Thm:Aux}, we finally obtain the desired gradient stability inequality as stated in Theorem \ref{M:S:THM}.

Finally, we consider $\Omega=\R^N$, $N \geq 2$, and investigate whether there exist positive constants $C$ and $\alpha$ such that for every $u\in W_0^{1,p}(\R^N),$  $$\frac{\big(
\kappa|\operatorname{supp}u|^{\frac1N-\frac1p}\|\nabla u\|_{L^p(\mathbb R^N)}\big)^p
-\|u\|_{L^\infty(\mathbb R^N)}^p}{
\|u\|_{L^\infty(\mathbb R^N)}^p
}\ge C \left(\inf_{a,b,x_0}
\frac{\left\|
\nabla u- \nabla v_{a,b,x_0}
\right\|_{L^p(\mathbb R^N)}}{\|  \nabla u\|_{L^p(\mathbb R^N)}}\right)^\alpha \,.$$
Unfortunately, this remains open. However, we could prove a pseudo-stability result as in the one-dimensional case (Theorem \ref{Morr-th-1d}), see Theorem \ref{Rn-stability}. Our proof is inspired by Cianchi's result in \cite{cian2008}, with standard modifications based on the reduction Lemma \ref{Thm:Aux:Rn}. For completeness, we include the details in the appendix (Theorem \ref{Rn-stability}).

\subsection{Hardy inequality} The second part of this article focuses on the quantitative stability of the $L^p$-Hardy inequality \eqref{Hardy-p}.  The non-trivial extremals of the $L^p$-Hardy inequality exist if $\Omega$ is a bounded $C^{1,\gamma}$-domain with a positive spectral gap $\Gamma_p(\Omega)$ (defined in Section \ref{sec-Hardy}), and in this case the space of extremals  is of the form
 $$ \mathcal{M}_{\mathrm{Hardy}}:=\{c\,\Psi: c\in \mathbb{R}\} \,,$$  where $\Psi>0$ solves the associated weighted eigenvalue problem \eqref{EVP} and $\|\nabla \Psi\|_{L^p(\Omega)}=1 $.

 As in the previous case, we define \textit{$p$-Hardy deficit} functional by  \[
\delta_{H}(u):=
\begin{cases}
\frac{\| \nabla u\|_{L^p(\Omega)}^{p'}- \mathbb{H}_{p}(\Omega)^{p'}\|u\|^{p'}_{L^p(\Omega,\,d_\Omega^{-p})}}{\|u\|^{p'}_{L^p(\Omega,\,d_\Omega^{-p})}}
& \text{if } p<2\,, \\[3ex]
\frac{\|\nabla u\|_{L^p(\Omega)}^{p}- \mathbb{H}_p(\Omega)^p\|u\|^{p}_{L^p(\Omega,\,d_\Omega^{-p})}}{\|u\|^{p}_{L^p(\Omega,\,d_\Omega^{-p})}}
& \text{if } p\ge 2 \,,
\end{cases}
\]
{where $\mathbb{H}_p(\Omega)$ is the  Hardy constant {for $C^{1,\gamma}$ bounded domains} (defined in Section \ref{sec-Hardy}).}
Define the distance of $u$ from the space of extremals $\mathcal{M}_{\mathrm{Hardy}}$ by 
\begin{equation}\label{dis_H}
  D_H(u)=\inf\left\{\frac{\|\nabla u-\nabla v\|_{L^p(\Omega)}}{\| u\|_{L^p(\Omega,\,d_\Omega^{-p})}}:v\in \mathcal{M}_{\mathrm{Hardy}}\right\},\quad u\in W^{1,p}_0(\Omega).  
\end{equation}

As we mentioned in the introduction, unlike the Morrey inequality, there exists $C>0$ such that $\delta_H(\cdot) \geq C D_H(\cdot)^2$ for the $L^2$-Hardy inequality. More precisely, we have the following. 
  \begin{thm}\label{Hardy-Thm-2} Let $p=2$ and $\Omega$ be a bounded  $C^{1,\gamma}$-domain in $\mathbb{R}^N$ such that $\Gamma_2(\Omega)>0.$ Then 
\begin{equation}\label{Hardy-2}
  \delta_H(u)\ge \left(1-\frac{\mathbb{H}_2(\Omega)^2}{\Lambda(\Omega)^2}\right) D_H(u)^2, \quad\forall\,u\in W^{1,2}_0(\Omega) \,,  
\end{equation}
where $$ \Lambda(\Omega):=\inf\left\{ \frac{\|\nabla u\|_{L^2(\Omega)}}{\|u\|_{L^2(\Omega,\,d_\Omega^{-p})}}: \, {u\in W^{1,2}_0(\Omega)}\setminus\{0\},\; \langle u,\Psi\rangle_{W^{1,2}_0(\Omega)}=0\right\} > \mathbb{H}_2(\Omega)  \,.$$
Moreover, the constant  $ \left(1-\frac{\mathbb{H}_2(\Omega)^2}{\Lambda(\Omega)^2}\right)$ is sharp. Furthermore, if $\Lambda(\Omega)$ is attained at $v,$ then the equality in \eqref{Hardy-2} occurs at $v+c\Psi$ for any $c \in \mathbb{R}$. Conversely,  if the equality in \eqref{Hardy-2} holds at $u$  then there exists  $c \in \mathbb{R}$ such that $\Lambda(\Omega)$ is attained at $u + c\Psi$.
\end{thm}

\begin{rem}
We provide a sufficient condition for the attainment of $\Lambda(\Omega)$. It turns out that $\Lambda(\Omega)$ is attained if the {\it Hardy constant at infinity} $\mathbb{H}_{2,\infty}(\Omega)> \Lambda(\Omega)$, see Proposition \ref{Prop-attain}.
\end{rem}
The proof of this theorem requires the standard Bianchi-Egnell-type arguments in stability theory. But, to show that $\Lambda(\Omega)> \mathbb{H}_2(\Omega)$ we use concentration compactness arguments. To prove the sharpness and the attainment of the constant  in \eqref{Hardy-2}, we use the Hilbert space structure of the underlying spaces and exploit the fact that the space of extremals is of one-dimension.

Finally, we treat the non-linear cases into two separate regimes $2<p<\infty$ and $1<p<2$ in Theorem \ref{Hardy-Thm-grad} and Theorem \ref{Hardy-Thm-grad<2}, respectively. Indeed, the non-linear cases are more involved. Moreover, to establish the quantitative gradient stability of the $L^p$-Hardy inequality in the two different regimes requires quite different {techniques}, as it can be seen explicitly in Section \ref{sec-Hardy}. 

\section{Stability of the Morrey inequality}\label{Morrey_section}
In this section, we will prove the quantitative pseduo-stability  (Theorem \ref{Morr-th-1d}) and the quantitative gradient stability (Theorem \ref{M:S:THM}) of the Morrey inequality.

Let us first recall the notion of Schwarz symmetrization.
\begin{defn}[Schwarz symmetrization]\label{defn-sym} Let $u$ be a measurable function on a domain $\Omega \subset \mathbb{R}^N$ such that the super-level sets of $|u|$ are of finite Lebesgue measure, i.e., $|\{|u|>t\}|<\infty$ for all $t>0$. Then the one-dimensional decreasing rearrangement $u^*:[0,|\Omega|) \rightarrow \mathbb{R}$ is defined as \[u^*(t)= \begin{cases}
\mathrm{ess.} \sup(|u|) \ \    \mbox{if} \ t=0\,,
    \\
    \inf \left\{s>0:|\{x\in \Omega:|u(x)|>s\}|<t\right\} \ \ \mbox{otherwise} \,,
    \end{cases}\] where $|E|$ denotes the Lebesgue measure of a measurable set $E\subset\mathbb{R}^N$.
    Let $\Omega_{\star}$ be the open ball in $\mathbb{R}^N$ centered at the origin and has the same Lebesgue measure as $\Omega$. The Schwarz symmetrization $u^\star: \Omega_{\star} \rightarrow \mathbb{R}$ is defined as the symmetric decreasing rearrangement as follows $$u^\star(x)=u^*(\omega_N|x|^N)\,.$$   
 \end{defn}
 It is not difficult to see that $u^\star$  is nonnegative, lower semi-continuous, and radially decreasing. Observe that the $L^{\infty}$-norm remains unchanged under the Schwarz symmetrization, i.e,  \begin{equation}\label{1}
\|u^\star\|_{L^{\infty}(\Omega_{\star})}=\|u^*\|_{L^{\infty}(0,|\Omega|)}=\|u\|_{L^{\infty}(\Omega)}.
\end{equation} In fact, it is well known that the $L^r$-norms of the functions are preserved under the Schwarz symmetrization, i.e., $\|u^\star\|_{L^r(\Omega)}=\|u^*\|_{L^r(0,|\Omega|)}=\|u\|_{L^r(\Omega)}$ for all $r \in [1,\infty]$. On the other hand, for $r \in (1,\infty)$, the P\'olya--Szeg\"o inequality says that the $L^r$-energy of the Sobolev-functions  gets reduced under the Schwarz symmetrization, i.e., 
\begin{equation}\label{polya}
    \int_{\Omega_{\star}}|\nabla u^\star|^r{\mathrm{d}}x\leq \int_{\Omega}|\nabla u|^r{\mathrm{d}}x\,, \quad \forall \, u\in  W^{1,r}_0(\Omega).
\end{equation} 

We will first prove the pseudo-stability of the Morrey inequality in a ball for the linear case, i.e., when $p=2$. Hence, $N=1$, i.e., in this case $\Omega$ is a symmetric interval centered at the origin. It is also known from \cite[Lemma 2.1]{Lindgren}  (see also \cite{Erco2016}) that $u\in\mathcal{M}_{\mathrm{Morrey}}$ if and only if
   \begin{equation}\label{weak-sol}
        \int_\Omega  u' v'\mathrm{d}x=\mathbb{M}_{2}(\Omega)^2u(0)v(0)\,,\quad\forall\,v\in  W^{1,2}_0(\Omega).
    \end{equation} 
\begin{proof}[Proof of Theorem \ref{Morr-th-1d}] 
Let $R>0$ and $\Omega=(-R,R).$
  WLOG, we can assume that
$\|u\|_{L^\infty(\Omega)}=\|u^\star\|_{L^\infty(\Omega)}=1.$ 
 Since $p=2$ and $N=1,$ $$\mathcal{M}_{\mathrm{Morrey}}=\left\{c\,\Phi: c \in \mathbb{R}\right\}, \ \ \mbox{where}\  \Phi=\left(R-|x|\right)/\sqrt{2R}.$$  Then $\|\Phi'\|_{L^2(\Omega)}=1$ and consequently $\Phi(0)=\frac{1}{\mathbb{M}_2(\Omega)}$. Recall that $$D_M(u)=\inf\left\{\|u'-c\,\Phi'\|_{L^2(\Omega)}:c\in\mathbb{R}\right\}.$$ 
Since $\mathcal{M}_{\mathrm{Morrey}}$ is a one-dimensional subspace of the Hilbert space $W^{1,2}_0(\Omega)$, it is closed. Thus, by the Hilbert projection theorem, $D_M(u)$ is uniquely attained \cite[Theorem 7.1.2, Remark 7.1.8]{kesh-f}, i.e., there exists a unique $a_0=a_0(u)\in\mathbb{R}$ such that $$D_M(u)=\| u'- a_0\Phi'\|_{L^2(\Omega)} \,.$$ Moreover, we can decompose $u$ as $$u=a_0\Phi+D_M(u)
\,w,$$ where  $\langle w,\Phi\rangle_{W^{1,2}_0(\Omega)}=0$ with $\| w'\|_{L^2(\Omega)}=1$. From \eqref{weak-sol} we get 
$$\mathbb{M}_2(\Omega)^2 w(0)\Phi(0)=\int_\Omega w' \Phi'\,{\mathrm{d}}x=0.$$ This gives $w(0)=0,$ and consequently,  $a_0=\frac{u(0)-D_M(u)\,w(0)}{\Phi(0)}=\mathbb{M}_2(\Omega)u(0).$
Now observe that
\begin{align*}
\int_\Omega| u'|^2{\mathrm{d}}x&=a_0^2\int_\Omega | \Phi'|^2{\mathrm{d}}x+2a_0D_M(u)\int_\Omega \Phi' w'{\mathrm{d}}x+D_M(u)^2\int_\Omega| w'|^2{\mathrm{d}}x    \\& =a_0^2+D_M(u)^2 =\mathbb{M}_2(\Omega)^2u(0)^2+D_M(u)^2.
\end{align*}
Therefore, \begin{align}\label{stability-equality}
 \delta_M(u)& = \|  u'\|_{L^2(\Omega)}^2-\mathbb{M}_2(\Omega)^2\|u\|^2_{L^\infty(\Omega)}=D_M(u)^2+\mathbb{M}_2(\Omega)^2 \Big(u(0)^2-\|u\|^2_{L^\infty(\Omega)}\Big).
\end{align}
Thus, for any function $u$ which has maxima at $0$, $$\delta_M(u)=D_M(u)^2.$$
Since $u^\star(0)=\|u^\star\|_{L^\infty(\Omega)}=\|u\|_{L^\infty(\Omega)}$,  $${\delta}_M(u^\star)= D_M(u^\star)^2.$$
By using P\'olya--Szeg\"o inequality \eqref{polya} on $\Omega=(-R,R)=\Omega_{\star}$ and \eqref{1}, we get $$\delta_M(u)=\|  u'\|_{L^2(\Omega)}^2-\mathbb{M}_2(\Omega)^2\|u\|^2_{L^\infty(\Omega)} \ge  \|  {u^{\star}}'\|_{L^2(\Omega)}^2-\mathbb{M}_2(\Omega)^2\|u^\star\|^2_{L^\infty(\Omega)} =  {\delta}_M(u^\star).$$
Hence, $$\delta_M(u)\ge  D_M(u^\star)^2.$$
This proves the pseudo-stability of the Morrey inequality.

\textbf{Sharpness and attainability.} 
Let $v\in\mathcal{M}_{\mathrm{Morrey}}\setminus\{0\}$ and $w\in W_0^{1,2}(\Omega)\setminus \mathcal{M}_{\mathrm{Morrey}}$ be such that 
\begin{equation}\label{assuption}
   w( 0)=0,\quad(v+w)^\star=(v+w), \quad\|v+w\|_{L^\infty(\Omega)}=\|v\|_{L^\infty(\Omega)}. 
\end{equation}
For example, take  $v(x)=R-|x|$ and $w(x)=|x|-(x^2/R)$ for $x \in \Omega=(-R,R)$. 
Now define $\phi:=v+w.$ Our aim is to prove that
 $\delta_M(\phi)=D_M(\phi^\star)^2.$ 
Since $w(0)=0,$ from \eqref{weak-sol} we get
$$\int_\Omega  v'w'{\mathrm{d}}x=\mathbb{M}_{2}(\Omega)^2v(0)w(0)=0.$$ This, together with \eqref{assuption} and $v \in\mathcal{M}_{\mathrm{Morrey}}$ gives
$$\int_\Omega|  \phi'|^2{\mathrm{d}}x-\mathbb{M}_2(\Omega)^2\|\phi\|^2_{L^\infty(\Omega)}=  \int_\Omega| v'|^2{\mathrm{d}}x+\int_\Omega| w'|^2{\mathrm{d}}x-\mathbb{M}_2(\Omega)^2\|v\|^2_{L^\infty(\Omega)}=  \int_\Omega| w'|^2{\mathrm{d}}x,$$
which implies that $\delta_M(\phi)=\left(\frac{\|w'\|_{L^2(\Omega)}}{\|\phi\|_{L^\infty(\Omega)}}\right)^2.$ On the other hand, \begin{align*}
    D_M(\phi^\star)^2=D_M(\phi)^2 &=\frac{1}{\|\phi\|_{L^\infty(\Omega)}^2}\inf_{c\in \mathbb{R}}\int_\Omega|\phi'-cv'|^2{\mathrm{d}}x =\frac{1}{\|\phi\|_{L^\infty(\Omega)}^2}\inf_{c\in \mathbb{R}}\int_\Omega|(1-c)v'+w'|^2{\mathrm{d}}x\\&=\frac{1}{\|\phi\|_{L^\infty(\Omega)}^2}\left(\inf_{c\in \mathbb{R}}(1-c)^2\int_\Omega|v'|^2{\mathrm{d}}x+\int_\Omega|w'|^2{\mathrm{d}}x\right)=\left(\frac{\|w'\|_{L^2(\Omega)}}{\|\phi\|_{L^\infty(\Omega)}}\right)^2\,.
\end{align*}
Hence, equality in the stability inequality  \eqref{morrey-2} is attained at $\phi,$
completing the proof. 
\end{proof}

\begin{rem} \label{Rem32} \rm  Observe that the equation  \eqref{stability-equality} implies that $$\delta_M(u) =D_M(u)$$
for all $ u \in W^{1,2}_0((-R,R))$ with $u(0)=\|u\|_{L^\infty((-R,R))}$. Indeed, we have seen in Remark \ref{Rem1} that the instability occurs due to the functions which attain their maximum outside the origin.
    
\end{rem}

Next, we consider the higher-dimensional non-linear case, i.e., $p>N\ge 1$. To prove Theorem \ref{M:S:THM}, we first establish an analogous stability inequality as in \eqref{simplify} with respect to the following distance 
\begin{equation}\label{eq:ci_dist}
   \beta(u):=\inf\left\{\frac{\|u-v\|_{L^\infty(\Omega)}}{\|u\|_{L^\infty(\Omega)}}: v\in\mathcal{M}_{\mathrm{Morrey}}, \,\|v\|_{L^\infty(\Omega)}=\|u\|_{L^\infty(\Omega)} \right\},\quad u\in W_0^{1,p}(\Omega). 
\end{equation}
Precisely, we prove the following:
 \begin{thm}\label{thm:MS:sup} Let $p>N\ge 1$ and $\Omega=B_R(x_0).$ Then there exists a positive constant $C=C(p,N,R)$  such that for every $u\in W^{1,p}_0(\Omega),$ 
\begin{align} \label{eq:thm_stab33}
    \frac{\|\nabla u\|_{L^p(\Omega)}^{p}- \mathbb{M}_p(\Omega)^p\|u\|^{p}_{L^\infty(\Omega)}}{\|u\|^{p}_{L^\infty(\Omega)}}\ge C
\begin{cases}
\beta(u)^p
& \text{if } p\ge 2, \\[2ex]
\min\{\beta(u)^p,\beta(u)^2\}
& \text{if } p<2.
\end{cases}
\end{align}
\end{thm}

Our proof of the above theorem requires a crucial {\it{capacity}} estimate for singleton sets inside the unit ball, which will be derived now.   For $p>N$ and $a\in \Omega:=B_1(0),$ the capacity of the set $\{a\}$ with respect to $\Omega$ is given by 
\begin{align}
    \tag{$\mathcal{P}$}\label{p1}
 \mathrm{Cap}_p(\{a\};\Omega):=\inf\left\{\int_{\Omega}|\nabla \phi|^p{\mathrm{d}}x:\,\phi\in W_0^{1,p}(\Omega),\;\phi(a)=1\right\}.
    \end{align}
In view of Remark \ref{Rem32}, it is relevant to consider the above minimization problem. Using standard variational arguments, one can show that the minimizers exist, say, $ u_a$ minimize \eqref{p1}, i.e., $$ \mathrm{Cap}_p(\{a\};\Omega)=\int_\Omega|\nabla  u_a|^p{\mathrm{d}}x. $$
We claim:
\begin{equation}
    \label{ms:e0} \mathrm{Cap}_p(\{a\};\Omega)\ge \mathbb{M}_p(\Omega)^p
\end{equation}
 and
\begin{equation}\label{ms:e3}
    \int_\Omega|\nabla u_a|^{p-2}\,\nabla u_a\cdot \nabla(u-u_a){\mathrm{d}}x=0\;\;\text{for every}\;u\in W_0^{1,p}(\Omega)\;\text{with}\;u(a)=1.
\end{equation} 
 Note that $\|u_a^{\star}\|_{L^\infty(\Omega)}=\|u_a\|_{L^\infty(\Omega)}\ge u_a(a)= 1$ (see \eqref{1}). Applying the Morrey inequality \eqref{morrey} for $u_a^{\star}$ and P\'olya--Szeg\"o inequality \eqref{polya}, we get $$\mathbb{M}_p(\Omega)^p\le\frac{1}{\|u_a^{\star}\|^p_{L^\infty(\Omega)}} \int_\Omega|\nabla u_a^{\star}|^p{\mathrm{d}}x\le  \int_\Omega|\nabla u_a^{\star}|^p{\mathrm{d}}x\le \int_\Omega|\nabla u_a|^p{\mathrm{d}}x=\mathrm{Cap}_p(\{a\};\Omega),$$ and conclude \eqref{ms:e0}.
Next, to prove \eqref{ms:e3}, let \(u\in W_0^{1,p}(\Omega)\) satisfy $u(a)=1$. Consider the functions $w_t:\Omega\to\mathbb{R}$ defined by
$$
w_t(x)=u_a(x)+t(u-u_a)(x),\qquad t\in\mathbb{R}.
$$
Since $w_t(a)=1$, each $w_t$ is admissible for \eqref{p1}, and hence $$\int_\Omega|\nabla u_a|^p{\mathrm{d}}x\le \int_\Omega|\nabla w_t|^p{\mathrm{d}}x\quad \forall\,t\in \R.$$  Thus, the function $F(t)=\int_\Omega|\nabla w_t|^p{\mathrm{d}}x$ has its minimum at $t=0.$ Consequently, $F'(0)=0,$ which proves our claim \eqref{ms:e3}.

The following capacity estimate improves \eqref{ms:e0}, which will be used in the proof of Theorem \ref{thm:MS:sup}. The estimate is also of independent interest.
\begin{prop} \label{Prop:cap}
    Let $p>N\ge 1$ and $a\in \Omega =B_1(0)$. Then we have   
    \begin{equation}\label{cap-estimate-p}
   \mathrm{Cap}_p(\{a\};\Omega)\ge \left(1-|a|^\frac{p-N}{p-1}\right)^{-p/p'} \mathbb{M}_p(\Omega)^p.
    \end{equation}
 Moreover, for $a=0$, equality holds in \eqref{cap-estimate-p}, namely,  $\mathrm{Cap}_p(\{0\};\Omega)=\mathbb{M}_p(\Omega)^p.$
\end{prop}
\begin{proof}  
  Let $\phi\in W^{1,p}_0(\Omega)$ be such that $\phi(a)=1.$ Then using polar decomposition $x=(r,\theta),$ $r=|x|$ and $\theta\in \mathbb{S}^{N-1},$ we obtain [cf. \cite[Lemma 2.1]{Drabek.}]
\begin{align*}
    1=\phi(a)& =\phi(|a|,\theta) =-\int_{|a|}^1\frac{\partial\phi}{\partial t}(t,\theta)\,\mathrm{d}t =-\int_{|a|}^1t^\frac{1-N}{p}t^\frac{N-1}{p}\frac{\partial\phi}{\partial t}(t,\theta)\,\mathrm{d}t\\&\le \left(\int_{|a|}^1t^\frac{(1-N)p'}{p}\,\mathrm{d}t\right)^{1/p'}\left(\int_{|a|}^1t^{N-1}\left|\frac{\partial\phi}{\partial t}(t,\theta)\right|^p\,\mathrm{d}t\right)^{1/p}\\&=\left(\frac{p-1}{p-N}\right)^{1/p'} \left(1-|a|^\frac{p-N}{p-1}\right)^{1/p'}\left(\int_{|a|}^1t^{N-1}\left|\frac{\partial\phi}{\partial t}(t,\theta)\right|^p\,\mathrm{d}t\right)^{1/p}\\&\le\left(\frac{p-1}{p-N}\right)^{1/p'} \left(1-|a|^\frac{p-N}{p-1}\right)^{1/p'}\left(\int_{0}^1t^{N-1}\left|\nabla\phi(t,\theta)\right|^p\,\mathrm{d}t\right)^{1/p},
\end{align*}
where we used H\"older inequality in the first inequality. This implies that $$\left(\frac{p-1}{p-N}\right)^{-p/p'} \left(1-|a|^\frac{p-N}{p-1}\right)^{-p/p'}\le \int_{0}^1t^{N-1}\left|\nabla\phi(t,\theta)\right|^p\,\mathrm{d}t.$$ Integrating both sides of the above inequality over $\mathbb{S}^{N-1}$ and thanks to the explicit value of the Morrrey constant $\mathbb{M}_p(\Omega)$, we obtain $$\mathbb{M}_p(\Omega)^p \left(1-|a|^\frac{p-N}{p-1}\right)^{-p/p'}=\left(\frac{p-1}{p-N}\right)^{-p/p'} \left(1-|a|^\frac{p-N}{p-1}\right)^{-p/p'}N\omega_N\le \int_\Omega|\nabla\phi|^p{\mathrm{d}}x.$$ 
Now, taking the infimum over all such $\phi,$ we obtain \eqref{cap-estimate-p}. Next, we prove that the equality holds for $a=0$. Since $\mathbb{M}_p(\Omega)$ is attained at $\Phi(x)=1-|x|^\frac{p-N}{p-1},$ $$\mathrm{Cap}_p(\{0\};\Omega)\le \int_\Omega|\nabla\Phi|^p{\mathrm{d}}x=\mathbb{M}_p(\Omega)^p.$$ Consequently, the equality of the above inequality follows from \eqref{ms:e0}. Hence the proof.
\end{proof}
We need the following elementary algebraic vector inequality 
to prove Theorem \ref{thm:MS:sup}. For any $1<p<\infty$ there exists a constant $C=C(p)>0$ such that
for all $X,Y\in\R^N$ we have
    \begin{equation}\label{p>2}
    |X+Y|^p\ge |X|^p+p|X|^{p-2}X\cdot Y+C
|Y|^p\quad
 \text{if } p\ge 2;
\end{equation}
\begin{equation}\label{p<2}
    |X+Y|^p\ge |X|^p+p|X|^{p-2}X\cdot Y+C
\min\{|Y|^p,|X|^{p-2}|Y|^2\}\quad
 \text{if } p<2.
 \end{equation}
For a proof, we refer to \cite[Lemma 3.1]{Tertikas}.

\begin{proof}[Proof of Theorem \ref{thm:MS:sup}] { For $1<p<\infty$, let us define
    $$\tilde{\delta}_M(u):=\frac{\|\nabla u\|_{L^p(\Omega)}^{p}- \mathbb{M}_p(\Omega)^p\|u\|^{p}_{L^\infty(\Omega)}}{\|u\|^{p}_{L^\infty(\Omega)}}\,,\quad u\in W^{1,p}_0(\Omega) \,,$$ i.e., $\tilde{\delta}_M(\cdot)$ is same as $\delta_M(\cdot)$ if $p\ge 2$, but for $p< 2$ they are not equal.} First of all, notice that it is enough to derive the desired inequality \eqref{eq:thm_stab33} for $\Omega=B_1(0)$, and without loss of generality we can take
 $$ \|u\|_{L^\infty(\Omega)}=1.$$
 Choose one of such $u$ arbitrarily. Since $W^{1,p}_0(\Omega)\hookrightarrow C^{0,\alpha}(\Omega),$ there exists $a\in \Omega:=B_1$ such that $u(a)=\|u\|_{L^\infty(\Omega)}=1.$ Let the capacity $\mathrm{Cap}_p(\{a\};\Omega)$ is attained at $u_a \in W^{1,p}_0(\Omega)$. 
We estimate the deficit $\tilde{\delta}_M(u)$ as
\begin{align}
        \label{ms:e4}
  \tilde{\delta}_M(u) &=\|\nabla u\|^p_{L^p(\Omega)}-\mathbb{M}_p(\Omega)^p  \nonumber \\ 
    & =[\|\nabla u\|^p_{L^p(\Omega)}-\|\nabla u_a\|^p_{L^p(\Omega)}]+[\mathrm{Cap}_p(\{a\};\Omega)-\mathbb{M}_p(\Omega)^p]
    =:I_a(u)+J(a).
\end{align} 
On the other hand, notice that 
\begin{equation}\label{ms:10}
    \beta(u)\le \|u-\Phi\|_{L^\infty(\Omega)}\,,\quad\text{where}\;\; \Phi(x):=1-|x|^\frac{p-N}{p-1}\in \mathcal{M}_{\mathrm{Morrey}}.
\end{equation} 
Our aim is to find suitable lower estimates of $I_a(u)$, $J(a)$, and an upper estimate of $\|u-\Phi\|_{L^\infty(\Omega)}$, which will be compatible in such a way that we would be able to establish \eqref{eq:thm_stab33}.
We will perform this job in several steps as follows.

\textbf{Step-1.}
First we estimate $I_a(u)$.

{Case-1.} $p\ge 2.$ Use \eqref{p>2} with $X=\nabla u_a$ and $Y=\nabla (u-u_a)$ to get 
$$|\nabla u|^p-|\nabla u_a|^p\ge p|\nabla u_a|^{p-2}\nabla u_a\cdot\nabla(u-u_a)+C|\nabla(u-u_a)|^p,$$ for some $C=C(p)>0.$ Integrate both sides of the above inequality over $\Omega$ and use \eqref{ms:e3} to get
\begin{equation}\label{ms:e5}
    I_a(u)\ge C\|\nabla(u-u_a)\|_{L^p(\Omega)}^p\ge C\,\mathbb{M}_p(\Omega)^p\|u-u_a\|^p_{L^\infty(\Omega)},
\end{equation}
 where we used the Morrey inequality \eqref{morrey} in the last inequality. 

 {Case-2.} $p<2$. In this case, the only admissible dimension is $N=1$. We will derive an analogue of \eqref{ms:e5} when $|a|<1/2$.  Use \eqref{p<2} with $X= u_a'$ and $Y= u'-u_a'$ to get 
$$| u'|^p-| u_a'|^p\ge p| u_a'|^{p-2} u_a'(u'-u_a')+c\,\min\{|u'-u_a'|^p, |u_a'|^{p-2}|u'-u_a'|^{2}\},$$ for some $c=c(p)>0.$ Integrate both sides of the above inequality over $\Omega$ and use \eqref{ms:e3} to obtain
\begin{align}\label{ms<2:e1}
    I_a(u) & \ge  c\int_\Omega \min\{|u'-u_a'|^p, |u_a'|^{p-2}|u'-u_a'|^{2}\}{\mathrm{d}}x  \nonumber\\
    & =  c \left[\int_A |u'-u_a'|^p{\mathrm{d}}x+\int_{A^c\cap\Omega}|u_a'|^{p-2}|u'-u_a'|^{2}{\mathrm{d}}x \right].
\end{align}
where $A=\{x\in \Omega: |u'-u_a'|>|u_a'|\}$.  Applying H\"older inequality, one gets 
\begin{equation}\label{ms:sup1}
\int_{A^c\cap\Omega} |u'-u_a'|^p{\mathrm{d}}x\le \left(\int_{A^c\cap\Omega}|u_a'|^{p-2}|u'-u_a'|^{2}{\mathrm{d}}x\right)^{p/2}\left(\int_{A^c\cap\Omega}|u_a'|^p{\mathrm{d}}x\right)^{(2-p)/2}\,,    
\end{equation}
 and observe that 
\begin{equation}\label{ms:sup2}
  \int_\Omega |u_a'|^p{\mathrm{d}}x=  \mathrm{Cap}_p(\{a\};\Omega)\le \|\zeta'\|_{L^p((-1,1))}^p=C_1,\quad \forall\, |a|<1/2
\end{equation}  where $\zeta\in C_c^\infty((-1,1))$ is a  function satisfying $\zeta\equiv 1$ in $[-\frac{1}{2},\frac{1}{2}]$.  Combining  \eqref{ms<2:e1}, \eqref{ms:sup1} and \eqref{ms:sup2}, we obtain $$I_a(u)\ge C_2\left[\int_{A} |u'-u_a'|^p{\mathrm{d}}x+\left(\int_{A^c\cap\Omega} |u'-u_a'|^p{\mathrm{d}}x\right)^{2/p} \right],$$ for some constant $C_2=C_2(p)>0.$ Using the elementary inequality $a+b^\frac{2}{p}\ge C_3\min\{(a+b),(a+b)^\frac{2}{p}\}$ for all $a,b\ge0,\,p\in (1,2),$ for some $C_3=C_3(p)>0$, we get
\begin{align*}
 \int_{A} |u'-u_a'|^p {\mathrm{d}}x+\left(\int_{A^c\cap\Omega} |u'-u_a'|^p{\mathrm{d}}x\right)^{2/p}\ge C_3\min\left\{\int_{\Omega}|u'-u_a'|^p{\mathrm{d}}x,\left( \int_{\Omega}|u'-u_a'|^p{\mathrm{d}}x\right)^{2/p}\right\}.
\end{align*} 
 Combining the above two inequalities and applying the Morrey inequality \eqref{morrey}, we obtain
 \begin{equation}\label{ms<2:e2}
     I_a(u)\ge C_4 \mathbb{M}_p(\Omega)^p\min\left\{\|u-u_a\|_{L^\infty(\Omega)}^p, \|u-u_a\|_{L^\infty(\Omega)}^2\right\}\,,\quad\forall\,|a|<1/2,
 \end{equation}
for some $C_4=C_4(p,N)>0.$

{\bf Step-2.} Estimate $J(a)$. By Proposition \ref{Prop:cap},   
$$J(a)=\mathrm{Cap}_p(\{a\};\Omega)-\mathbb{M}_p(\Omega)^p \ge \mathbb{M}_p(\Omega)^p\left[\left(1-|a|^\frac{p-N}{p-1}\right)^{-p/p'}-1\right].$$ 
Applying the Bernoulli's inequality: $$(1-t)^{-p/p'}=(1-t)^{-(p-1)}\ge 1+(p-1)t,\quad t\in (0,1)$$ with $t=|a|^\frac{p-N}{p-1}$ in the above estimate of $J(a)$, we get 
\begin{equation}\label{ms:e6}
    J(a)\ge (p-1)\mathbb{M}_p(\Omega)^p\,|a|^\frac{p-N}{p-1}\,,\qquad \forall\,p>1,\,|a|<1.
\end{equation}

\textbf{Step-3.}  Estimate $\|u-\Phi\|_{L^\infty(\Omega)}$. By the triangle inequality 
$$\|u-\Phi\|_{L^\infty(\Omega)}\le \|u-u_a\|_{L^\infty(\Omega)}+\|u_a-\Phi\|_{L^\infty(\Omega)}.$$
We already have an upper estimate of $\|u-u_a\|^p_{L^\infty(\Omega)}$ by $I_a(u)$ due to \eqref{ms:e5}. Now, we estimate $\|u_a-\Phi\|_{L^\infty(\Omega)}$. Consider $$w_a(x):=\frac{\Phi(x)}{\Phi(a)},\quad x\in \Omega\,,$$ where $\Phi$ is given by \eqref{ms:10}. Then the triangle inequality gives  $$\|u_a-\Phi\|_{L^\infty(\Omega)}\le \|u_a-w_a\|_{L^\infty(\Omega)}+\|w_a-\Phi\|_{L^\infty(\Omega)}.$$
The motivation for considering the function $w_a$ is that it is explicitly known and satisfies $w_a(a)=1$, which helps us to obtain the desired estimates of $\|u_a-w_a\|_{L^\infty(\Omega)}$ and $\|w_a-\Phi\|_{L^\infty(\Omega)}.$ Note that in Step-1, we estimate $\|u-u_a\|_{L^\infty(\Omega)}$ (see \eqref{ms:e5}, \eqref{ms<2:e2}) by using \eqref{ms:e3}, which is proved by using the fact that $u(a)=1$. Thus, in a similar way, we can prove that \eqref{ms:e5} holds with $w_a$ in place of $u$, i.e.,
\begin{equation}\label{ms:e7}
    \|u_a-w_a\|^p_{L^\infty(\Omega)}\le \frac{1}{C\mathbb{M}_p(\Omega)^p}[\|\nabla w_a\|^p_{L^p(\Omega)}-\|\nabla u_a\|^p_{L^p(\Omega)}],\quad \forall\,p\ge 2.
\end{equation}
 Similarly, we can also prove that  \eqref{ms<2:e2} holds with $w_a$ in place of $u$, i.e., for $|a|<1/2$
\begin{equation}\label{ms<2:e3}
    \min\{\|u_a-w_a\|^p_{L^\infty(\Omega)}, \|u_a-w_a\|^2_{L^\infty(\Omega)}\}\le \frac{1}{C_4\mathbb{M}_p(\Omega)^p}[\|\nabla w_a\|^p_{L^p(\Omega)}-\|\nabla u_a\|^p_{L^p(\Omega)}],\;\;  \forall\,p<2 \,.
\end{equation}
To estimate the right hand side of \eqref{ms:e7} and \eqref{ms<2:e3}, observe that $$\|\nabla w_a\|^p_{L^p(\Omega)}=\frac{\|\nabla \Phi\|^p_{L^p(\Omega)}}{\Phi(a)^p}=\frac{\mathbb{M}_p(\Omega)^p}{(1-|a|^\frac{p-N}{p-1})^p} \quad \mbox{and}  \quad \|\nabla u_a\|^p_{L^p(\Omega)}=\mathrm{Cap}_p(\{a\};\Omega)\ge \mathbb{M}_p(\Omega)^p \,,$$  where the second estimate is due to \eqref{ms:e0}. Thus, we obtain 
$$\|\nabla w_a\|^p_{L^p(\Omega)}-\|\nabla u_a\|^p_{L^p(\Omega)}\le \mathbb{M}_p(\Omega)^p\left(\frac{1}{(1-|a|^\frac{p-N}{p-1})^p}-1\right) \le p\mathbb{M}_p(\Omega)^p\frac{|a|^\frac{p-N}{p-1}}{(1-|a|^\frac{p-N}{p-1})^p}$$
by using the mean value theorem in the last inequality. 
Consequently, from \eqref{ms:e7} we get 
\begin{equation}\label{ms:e8}
    \|u_a-w_a\|^p_{L^\infty(\Omega)}\le \frac{p|a|^\frac{p-N}{p-1}}{C(1-|a|^\frac{p-N}{p-1})^p}\,,\quad \forall\,p\ge 2.  
\end{equation}
and from \eqref{ms<2:e3} we get for $|a|<1/2$ \begin{equation}\label{ms<2:e4}
  \min\{  \|u_a-w_a\|^p_{L^\infty(\Omega)},  \|u_a-w_a\|^2_{L^\infty(\Omega)}\}\le \frac{p|a|^\frac{p-N}{p-1}}{C_4(1-|a|^\frac{p-N}{p-1})^p}\,,\quad \forall\,p<2.  
\end{equation}
Finally, we estimate $\|w_a-\Phi\|_{L^\infty(\Omega)}$ below: 
\begin{equation}\label{ms:e9}
    \|w_a-\Phi\|_{L^\infty(\Omega)}=\left\|\frac{1-\Phi(a)}{\Phi(a)}\Phi\right\|_{L^\infty(\Omega)}=\frac{|a|^\frac{p-N}{p-1}}{1-|a|^\frac{p-N}{p-1}}\,,\quad \,p>1.
\end{equation}

\textbf{Step-4.} Now we have all the required estimates for establishing the desired stability inequality \eqref{eq:thm_stab33}.  
Combining \eqref{ms:e4}, \eqref{ms:e5}, and \eqref{ms:e6} we get the lower estimate of $\tilde{\delta}_M(u)$:
\begin{equation}\label{ms:11}
   \tilde{\delta}_M(u)\ge \mathbb{M}_p(\Omega)^p(C+p-1)\,[\|u-u_a\|^p_{L^\infty(\Omega)}+|a|^\frac{p-N}{p-1}],\quad p\ge 2.
\end{equation}
Similarly, combining \eqref{ms:e4}, \eqref{ms<2:e2}, and \eqref{ms:e6}, we get the lower estimate of $\tilde{\delta}_M(u)$ for $|a|<1/2$:
\begin{equation}\label{ms<2:11}
   \tilde{\delta}_M(u)\ge \mathbb{M}_p(\Omega)^p(C_4+p-1)\,\big[\min\{\|u-u_a\|^p_{L^\infty(\Omega)},\|u-u_a\|^2_{L^\infty(\Omega)}\}+|a|^\frac{p-N}{p-1}\big],\quad p<2.
\end{equation}  Now note that $\|w_a-\Phi\|_{L^\infty(\Omega)} \to \infty$ if $|a|\to 1$ (see \eqref{ms:e9}).
To overcome this difficulty, we consider two cases. 
 
Case 1: $1/2\le |a|\le 1.$ Combine \eqref{ms:e4} and \eqref{ms:e6} to get $$ \tilde{\delta}_M(u)\ge J_a(u)\ge (p-1)\mathbb{M}_p(\Omega)^p |a|^\frac{p-N}{p-1}\ge (p-1)\mathbb{M}_p(\Omega)^p (1/2)^\frac{p-N}{p-1}.$$ On the other hand, \eqref{ms:10} gives  $$\beta(u)\le \|u-\Phi\|_{L^\infty(\Omega)} \le \|u\|_{L^\infty(\Omega)}+\|\Phi\|_{L^\infty(\Omega)}=2,$$ and consequently, the stability inequality \eqref{eq:thm_stab33} holds.

Case 2: $|a|<1/2.$ First, let us consider the case $p\ge 2$. 
One can see that
$$\|u-\Phi\|_{L^\infty(\Omega)}^p\le C_5( \|u-u_a\|_{L^\infty(\Omega)}^p+\|u_a-w_a\|_{L^\infty(\Omega)}^p+\|w_a-\Phi\|_{L^\infty(\Omega)}^p)$$
for some $C_5=C_5(p)>0$. Now, using \eqref{ms:e8} and \eqref{ms:e9}, we get $C_6=C_6(p,N)>0$ such that $$\|u-\Phi\|_{L^\infty(\Omega)}^p\le C_6\left(\|u-u_a\|_{L^\infty(\Omega)}^p+|a|^\frac{p-N}{p-1}+|a|^\frac{p(p-N)}{p-1}\right)\le 2C_6\left(\|u-u_a\|_{L^\infty(\Omega)}^p+|a|^\frac{p-N}{p-1}\right).$$
 Combining \eqref{ms:10}, the above inequality, and \eqref{ms:11} we obtain $$\beta(u)^p\le \|u-\Phi\|_{L^\infty(\Omega)}^p\le  C_7\tilde{\delta}_M(u),$$  for some $C_7=C_7(p,N)>0,$ i.e., \eqref{eq:thm_stab33} holds.
 Next, we consider the case $p<2$. Since $$\|u-\Phi\|_{L^\infty(\Omega)}\le  \|u-u_a\|_{L^\infty(\Omega)}+\|u_a-w_a\|_{L^\infty(\Omega)}+\|w_a-\Phi\|_{L^\infty(\Omega)},$$  one can verify that there exists $C_8=C_8(p)>0$ such that 
\begin{multline}\label{ms<2:e5}
\min\{\|u-\Phi\|_{L^\infty(\Omega)}^p, \|u-\Phi\|_{L^\infty(\Omega)}^2\}\le C_8\big[ \min\{\|u-u_a\|_{L^\infty(\Omega)}^p, \|u-u_a\|_{L^\infty(\Omega)}^2\}+\\\min\{\|u_a-w_a\|_{L^\infty(\Omega)}^p, \|u_a-w_a\|_{L^\infty(\Omega)}^2\}+\|w_a-\Phi\|_{L^\infty(\Omega)}^p\big].   
\end{multline}
Moreover, use \eqref{ms<2:e4} and \eqref{ms:e9} to obtain  \begin{equation}\label{ms<2:e15}
    \min\{\|u_a-w_a\|_{L^\infty(\Omega)}^p, \|u_a-w_a\|_{L^\infty(\Omega)}^2\}+\|w_a-\Phi\|_{L^\infty(\Omega)}^p\le C_9(|a|^\frac{p-N}{p-1}+|a|^\frac{p(p-N)}{p-1})\le 2C_9|a|^\frac{p-N}{p-1},
\end{equation}
  for some $C_9=C_9(p,N)>0.$ Now, combining \eqref{ms:10}, \eqref{ms<2:e5}, \eqref{ms<2:e15}, and \eqref{ms<2:11} we obtain $$\min\{\beta(u)^p,\beta(u)^2\}\le \min\{\|u-\Phi\|_{L^\infty(\Omega)}^p, \|u-\Phi\|_{L^\infty(\Omega)}^2\}\le C_{10}\,\tilde{\delta}_M(u),$$ for some $C_{10}=C_{10}(p,N)>0,$ where we used the fact that $\min\{a^p,a^2\}\le \min\{b^p,b^2\}$ for all $b\ge a\ge 0$ in the first inequality. That is, \eqref{eq:thm_stab33} holds for $p<2$. Thus, the combination of Case-1 and Case-2 completes the proof.
\end{proof}

Next, we establish a reduction lemma that, together with Theorem \ref{thm:MS:sup}, yields Theorem \ref{M:S:THM}.  

\begin{lem}\label{Thm:Aux}
Let $p>N\ge 1$ and  $\Omega \subset $ be a ball in $\mathbb{R}^N$. Then 
\begin{equation}\label{eqn-p-case}
  \left(  \frac{\|\nabla u-\nabla v\|_{L^p(\Omega)}}{\|u\|_{L^\infty(\Omega)}}\right)^\alpha\le C_1\delta_M(u)+C_2 \frac{\| u- v\|_{L^\infty(\Omega)}}{\|u\|_{L^\infty(\Omega)}}
\end{equation}
for all $u \in W^{1,p}_0(\Omega)$ and for any $v \in \mathcal{M}_{\mathrm{Morrey}}$ with $\|u\|_{L^\infty(\Omega)} = \|v\|_{L^\infty(\Omega)}$,
where $\alpha=p'$ when $p<2,$ and  $\alpha=p$ when $p\ge 2.$ Here $C_1=2^\alpha$ and $C_2=\alpha 2^{\alpha-1}\mathbb{M}_p(\Omega)^\alpha.$
\end{lem}
\begin{proof}
We mainly follow the ideas of the proof of Theorem 1.1 in \cite{Neu2020} to prove \eqref{eqn-p-case}. We will use 
the following version of Clarkson’s inequalities for vector-valued functions: for all $F,G:\mathbb{R}^N\rightarrow\mathbb{R}^N$ with $|F|,|G|\in L^p(\mathbb{R}^N)$ we have 
    \begin{equation}\label{1-2}
        \left\|\frac{F+G}{2}\right\|_{L^{p}(\mathbb{R}^N)}^{p'}+\left\|\frac{F-G}{2}\right\|_{L^{p}(\mathbb{R}^N)}^{p'}\le \left(\frac{\|F\|_{L^{p}(\mathbb{R}^N)}^p+\|G\|_{L^{p}(\mathbb{R}^N)}^p}{2}\right)^{p'/p},
    \end{equation}
    if $p<2,$ and  \begin{equation}\label{p-2}
        \left\|\frac{F+G}{2}\right\|_{L^{p}(\mathbb{R}^N)}^{p}+\left\|\frac{F-G}{2}\right\|_{L^{p}(\mathbb{R}^N)}^{p}\le \frac{\|F\|_{L^{p}(\mathbb{R}^N)}^p+\|G\|_{L^{p}(\mathbb{R}^N)}^p}{2},
    \end{equation}
if $p\ge 2.$ For a proof of \eqref{1-2} and \eqref{p-2}, see \cite[Section 3]{Neu2020}.

We first address the case $p \in (1, 2)$. Applying \eqref{1-2} with $F = \nabla u$  and
$G = \nabla v$, after extending both functions by zero outside $\Omega$, we deduce that
    $$ \left\|\frac{\nabla u-\nabla v}{2}\right\|_{L^p(\Omega)}^{p'}\le \left(\frac{\|\nabla u\|_{L^p(\Omega)}^p+\|\nabla v\|_{L^p(\Omega)}^p}{2}\right)^{p'/p}-\left\|\frac{\nabla u+\nabla v}{2}\right\|_{L^p(\Omega)}^{p'}.$$ 
By the Morrey inequality \eqref{morrey}, we further obtain $$\| \nabla u\|_{L^p(\Omega)}\ge \mathbb{M}_p(\Omega)\|u\|_{L^\infty(\Omega)}=\mathbb{M}_p(\Omega)\|v\|_{L^\infty(\Omega)}=\| \nabla v\|_{L^p(\Omega)},$$
and $$\|\nabla u+\nabla v\|_{L^p(\Omega)}\ge\mathbb{M}_p(\Omega) \|u+v\|_{L^\infty(\Omega)}.$$ Combining the above three inequalities, we obtain
\begin{align}\label{1-thmp}
   \left\|\frac{\nabla u-\nabla v}{2}\right\|_{L^p(\Omega)}^{p'}&\le \|\nabla u\|_{L^p(\Omega)}^{p'}-\mathbb{M}_p(\Omega)^{p'}\left\|\frac{ u+v}{2}\right\|_{L^\infty(\Omega)}^{p'}.
\end{align}
  Next, we show that 
  \begin{equation}\label{2-thmp}
      \left\|\frac{ u+v}{2}\right\|_{L^\infty(\Omega)}^{p'}\ge \|u\|_{L^\infty(\Omega)}^{p'}-p'\|u\|_{L^\infty(\Omega)}^{p'-1}\left\|\frac{ u-v}{2}\right\|_{L^\infty(\Omega)}.
  \end{equation}
Use Minkowski's inequality to obtain
$$ \left\|\frac{ u+v}{2}\right\|_{L^\infty(\Omega)}^{p'}\ge  \left(\|u\|_{L^\infty(\Omega)}-\left\|\frac{ u-v}{2}\right\|_{L^\infty(\Omega)}\right)^{p'}.$$
Using the convexity of the function $t\mapsto t^{p'}$, we obtain 
  $$\left(\|u\|_{L^\infty(\Omega)}-\left\|\frac{ u-v}{2}\right\|_{L^\infty(\Omega)}\right)^{p'}\ge \|u\|_{L^\infty(\Omega)}^{p'}-p'\|u\|_{L^\infty(\Omega)}^{p'-1}\left\|\frac{ u-v}{2}\right\|_{L^\infty(\Omega)}.$$ Clearly, \eqref{2-thmp} follows by combining the last two inequalities.  Finally, combining \eqref{1-thmp} and \eqref{2-thmp} and dividing 
by $\|u\|_{L^\infty(\Omega)}^{p'}$ we get
$$\left(\frac{\|\nabla u-\nabla v\|_{L^p(\Omega)}}{2\|u\|_{L^\infty(\Omega)}}\right)^{p'}\le\frac{ \|\nabla u\|_{L^p(\Omega)}^{p'}-\mathbb{M}_p(\Omega)^{p'}\|u\|_{L^\infty(\Omega)}^{p'}}{{\|u\|^{p'}_{L^\infty(\Omega)}}}+p'\mathbb{M}_p(\Omega)^{p'}\frac{ \|u-v\|_{L^\infty(\Omega)}}{2\|u\|_{L^\infty(\Omega)}} \,,$$ which establishes
 \eqref{eqn-p-case} with $C_1=2^{p'}$ and $C_2=p'2^{p'-1}\mathbb{M}_p(\Omega)^{p'}$.

 The proof for case $ p \ge  2$ is completely analogous, and hence omitted. 
\end{proof}
Now, we are ready to prove Theorem \ref{M:S:THM}. The proof follows by applying stability with distance function $\beta(\cdot)$ (Theorem \ref{thm:MS:sup}) and the reduction lemma, Lemma \ref{Thm:Aux}. 

\begin{proof}[Proof of Theorem \ref{M:S:THM}] Let $\Omega=B_R(x_0),$ $u\in W^{1,p}_0(\Omega)$ and $\beta(u)$ be as in \eqref{eq:ci_dist}.
If $\beta(u) =0$, then it follows that $u \in \mathcal{M}_{\mathrm{Morrey}}$. Consequently,  $D_M(u)=0.$ Hence the theorem follows trivially. Now assume that $\beta(u) \ne 0$. Let $v\in \mathcal{M}_{\mathrm{Morrey}}$ be such that $$\|u\|_{L^\infty(\Omega)} = \|v\|_{L^\infty(\Omega)}\quad\text{and}\quad\frac{\| u- v\|_{L^\infty(\Omega)}}{\|u\|_{L^\infty(\Omega)}}\le 2\beta(u).$$  
Then letting $\alpha=p'$ for $p\in (1,2)$ and $\alpha=p$
 for $p\ge 2,$ we apply  Lemma \ref{Thm:Aux} to obtain
 \begin{align}\label{M:S:e1}
   D_M(u)^\alpha\le\left(\frac{\|\nabla u-\nabla v\|_{L^p(\Omega)}}{\|u\|_{L^\infty(\Omega)}}\right)^\alpha\le C\delta_M(u)+C\frac{\| u- v\|_{L^\infty(\Omega)}}{\|u\|_{L^\infty(\Omega)}}\le C\delta_M(u)+2C\beta(u),
 \end{align}
for some constant $C>0$ depending only on $p,\, N,$ and $R$

\textbf{Case 1.}  $p\ge 2.$
 By Theorem \ref{thm:MS:sup} we get a constant $C_1=C_1(p,N,R)>0$ such that $\delta_M(u)\ge C_1\beta(u)^p.$ Consequently from \eqref{M:S:e1} we obtain
\begin{equation*}
 D_M(u)^p\le C\delta_M(u)+2CC_1^{-1/p} \delta_M(u)^\frac{1}{p}\le C_2\max\{\delta_M(u),\delta_M(u)^\frac{1}{p}\},
\end{equation*}
for some constant $C_2=C_2(p,N,R)>1.$ Now using the fact that $a^p\le C_2\max\{b,b^{1/p}\} , \;a,b\ge 0$ implies $b\ge C_2^{-p}\min\{a^p,a^{p^2}\},$ we obtain 
$$ \delta_M(u)\ge C_2^{-p}\min\{D_M(u)^p,D_M(u)^{p^2}\}.$$ 

 \textbf{Case 2.}  $p\in (1,2).$ First, we verify that 
\begin{equation}\label{Ms:m<2}
  \frac{\|\nabla u\|_{L^p(\Omega)}^p-\mathbb{M}_p(\Omega)^p\|u\|^p_{L^\infty(\Omega)}}{\|u\|^p_{L^\infty(\Omega)}} \le \mathbb{M}_p(\Omega)^{p-p'}\,\delta_M(u).
\end{equation}
To see this, note that for any $p\in ( 1,2),$ the function $t\mapsto t^{p'}-t^p\quad\text{for} \, t\ge 1$ is an increasing function. So, if $a\ge b\ge1,$ we have 
\begin{equation}\label{MsEq2;m<2}
 a^p-b^p \leq   a^{p'}-b^{p'}  \,.
\end{equation}
Taking $a=\frac{\|\nabla u\|_{L^p(\Omega)}}{\mathbb{M}_p(\Omega)\|u\|_{L^\infty(\Omega)}}$ and $b=1$ in \eqref{MsEq2;m<2}, we obtain \eqref{Ms:m<2}.
Thus, by Theorem \ref{thm:MS:sup} we obtain $ \delta_M(u)\ge c  \min\left\{\beta(u)^p,\beta(u)^2\right\}$ for some $c=c(p,N,R)>0,$  
which implies that
\begin{equation}\label{MS:eq3<2}
   \beta(u)\le C_3 \max\left\{\delta_M(u)^\frac{1}{p},\delta_M(u)^\frac{1}{2}\right\},
\end{equation}
for some $C_3=C_3(p,N,R)>0.$ Combining \eqref{M:S:e1} and \eqref{MS:eq3<2}, we obtain
 \begin{align*}
   D_M(u)^{p'}\le C\delta_M(u)+2CC_3\max\left\{\delta_M(u)^\frac{1}{p},\delta_M(u)^\frac{1}{2}\right\}\le C_4\max\left\{\delta_M(u),\delta_M(u)^\frac{1}{2}\right\}, 
 \end{align*}
for some $C_4=C_4(p,N,R)>1$. This gives $$ \delta_M(u)\ge C_4^{-p'}\min\{D_M(u)^{p'},D_M(u)^{2p'}\}.$$ Hence the stability inequality \eqref{simplify} follows by combining Case 1 and Case 2.
\end{proof}

\section{Stability of the Hardy inequality}\label{sec-Hardy} In this section, we study the quantitative stability of the geometric $L^p$-Hardy inequality \eqref{Hardy-p} in a $C^{1,\gamma}$ bounded domain $\Omega \subset \mathbb{R}^N$. 
The $L^p$-Hardy constant of $\Omega$ 
will be
denoted by  
\begin{equation*}
    \mathbb{H}_p(\Omega):=\inf_{u\in W^{1,p}_0(\Omega)\setminus\{0\}} R_p(u)\,,  \  \ \ \ \mbox{where} \ R_p(u):=\frac{\|\nabla u\|_{L^p(\Omega)}}{\|u\|_{L^p(\Omega,\,d_\Omega^{-p})}} \,.
\end{equation*} 
If $\mathbb{H}_p(\Omega)$ is attained at some $\Psi\in W^{1,p}_0(\Omega)\setminus\{0\}$, then $\Psi$ is an extremal for the corresponding $L^p$-Hardy inequality in $\Omega$. Moreover, such an extremal (if exists) satisfies the corresponding Euler-Lagrange equation
\begin{align} \label{EVP}
    -\Delta_p u = \mathbb{H}_p(\Omega)^p\frac{|u|^{p-2}u}{d_\Omega^p} \quad\text{in}\,\ \Omega\,, \ \ \ \ u=0 \ \ \mbox{on} \  \partial \Omega \,.
\end{align} 

Our goal is to investigate the gradient stability of the $L^p$-Hardy inequalities in $\Omega$ when $\mathbb{H}_p(\Omega) $ is attained. Of course, there are domains for which $\mathbb{H}_p(\Omega) $ is not attained, i.e., $\mathcal{M}_{\mathrm{Hardy}}=\{0\}$.  
Let us recall a necessary and sufficient condition for the existence of non-trivial extremals for the $L^p$-Hardy inequality in bounded  $C^{1,\gamma}$-domains, where $\gamma\in (0,1].$ For this, we consider the \textit{Hardy constant at infinity} 
\begin{align*}
\mathbb{H}_{p,\infty}(\Omega)
=
\sup_{K\Subset\Omega} \inf \Bigl\{ R_p(u) : u\in W^{1,p}_0(\Omega\setminus \overline{K}) \setminus \{0\}
\Bigr\}.
\end{align*}
We write $A\Subset O $ if $\overline{A}$ is compact and $\overline{A}\subset O.$ For bounded $C^{1,\gamma}$-domains, $0<\mathbb{H}_{p}(\Omega) \leq \mathbb{H}_{p,\infty}(\Omega) = \frac{p-1}{p}$ \cite{Pinch19}.
It is known that \cite[Theorem 4.1 \& 5.1]{Pinch19} (see \cite{Pin-Ujjal} for the weighted case) $\mathbb{H}_p(\Omega)$ is attained if and only if the \textit{spectral gap} $$\Gamma_p(\Omega):=\mathbb{H}_{p,\infty}(\Omega)-\mathbb{H}_p(\Omega) >0 \,.$$
Furthermore, it is well known that if  $\mathbb{H}_p(\Omega)$ is attained then $\mathbb{H}_p(\Omega)^p$ is the first (i.e., the least) weighted eigenvalue of the problem \eqref{EVP}, and it is also principal (i.e., the corresponding eigenfunctions do not change their sign). In fact, the corresponding eigenfunctions are of definite sign (i.e., either strictly positive or negative). Moreover, the weighted eigenvalue $\mathbb{H}_p(\Omega)^p$ is simple (i.e.,  the corresponding eigenfunctions are constant multiples of each other).

The next proposition shows that we cannot expect a gradient stability inequality for the $L^p$-Hardy inequality if the spectral gap is zero.

\begin{prop}\label{Rem-Hardy} Let $\Omega$ be a $C^{1,\gamma}$ bounded domain with $\Gamma_p(\Omega)=0$. Then there do not exist $C_1,\alpha_1>0$ such that  $$\delta_H(u)\ge C_1 D_H(u)^{\alpha_1}\quad\forall\,u\in W_0^{1,p}(\Omega)\,,$$  i.e.,  gradient stability fails. Furthermore, if $\Omega=B_1(0)$ then there do not exist $C_2,\alpha_2>0$ such that  $$ \delta_H(u)\ge C_2 D_H(u^\star)^{\alpha_2}\quad\forall\,u\in W_0^{1,p}(\Omega)\,,$$  i.e.,  pseudo-gradient stability fails.
\end{prop}
\begin{proof} In the first case, $\mathcal{M}_{\mathrm{Hardy}}=\{0\}$ (see \cite[Theorem 5.1]{Pinch19}). Let $(u_k)$ be a minimizing sequence for $\mathbb{H}_p(\Omega).$ Then we have $$D_H(u_k)=\frac{\|\nabla u_k\|_{L^p(\Omega)}}{\|u_k\|_{L^p(\Omega,\,d_\Omega^{-p})}}\to \mathbb{H}_p(\Omega),\qquad \delta_H(u_k)\to 0\;\;\text{as}\;k\to 0.$$ Consequently, the gradient stability fails. 

For the second case, $\mathcal{M}_{\mathrm{Hardy}}=\{0\}$ 
as $\mathbb{H}_p(\Omega)=\frac{p-1}{p}=\mathbb{H}_{p,\infty}(\Omega)$ for $\Omega=B_1(0)$. 
For $\epsilon\in (0,1),$ consider
$$u_\epsilon(x):=(1-|x|)^{\frac{p-1}{p}+\epsilon}\quad \mbox{ in } \Omega.$$  Then  $\|\nabla u_{\epsilon}\|_{L^p(\Omega)}^p=\left(\frac{p-1}{p}+\epsilon\right)^pI_\epsilon $ and $\|u_\epsilon\|_{L^p(\Omega,\,d_\Omega^{-p})}^p=I_\epsilon,$ where $I_\epsilon=\int_\Omega(1-|x|)^{-1+p\epsilon}\,{\mathrm{d}}x.$ Thus, $$\delta_H(u_\epsilon)\le C \epsilon\; \ \;\text{and} \ \;\;D_H(u_{\epsilon}^\star)=D_H(u_\epsilon)=\frac{\|\nabla u_\epsilon\|_{L^p(\Omega)}}{\|u_\epsilon\|_{L^p(\Omega,\,d_\Omega^{-p})}}=\frac{p-1}{p}+\epsilon\,,\quad C=C(p)>0.$$ 
Therefore, for any $\alpha>0,$   $\frac{\delta_H(u_\epsilon)}{D_H(u_{\epsilon}^\star)^\alpha}\to 0$ taking $\epsilon\to 0$. Hence, pseudo-gradient stability fails.\end{proof}

In view of the above proposition,
{\bf throughout this section, we will assume that 
$\Omega$ is a bounded  $C^{1,\gamma}$ domain for some $\gamma \in (0,1]$ with $\Gamma_p(\Omega)>0$}. 

From the discussion before the proposition, we see that the space of extremals of the $L^p$-Hardy inequality in such $\Omega$ will be of the form
 $$ \mathcal{M}_{\mathrm{Hardy}}:=\{c\,\Psi: c\in \mathbb{R}\} \,,$$  where $\Psi>0$ solves \eqref{EVP} and $\|\nabla \Psi\|_{L^p(\Omega)}=1 $.

\subsection{The linear case: $p=2$}
First, we study the quantitative stability of the $L^2$-Hardy inequality with respect to the distance function $D_H(\cdot)$.
Recall that $$ \Lambda(\Omega)=\inf\left\{ R_2(u): \, u\in W^{1,2}_0(\Omega),\; \langle u,\Psi\rangle_{W^{1,2}_0(\Omega)}=0\right\} \,.$$ Before going to the proof of Theorem \ref{Hardy-Thm-2}, let us first show that $\Lambda(\Omega)>\mathbb{H}_2(\Omega)$, which will ensure that Theorem \ref{Hardy-Thm-2} indeed gives a stability result. Towards this, we use the concentration compactness lemma proved by Smets \cite{smets} (see also \cite{DY,Hou} for a variant of it).
Consider the following quantities 
$$S_{r}(x):
=
\inf\left\{ \|\nabla u\|^2_{L^2(\Omega)}:
u\in W_0^{1,2}(\Omega\cap B_r(x)),
\ \|u\|_{L^2(\Omega,\,d_{\Omega}^{-2})}=1
\right\},$$
 $$ S(x)  := \sup_{r>0} S_{r}(x)= \lim_{r\to 0} S_{r}(x)
\ \ \ \ \mbox{and} \ \ \ \ 
 S^*  :=  \inf_{x \in \overline{\Omega}} S(x) \,.$$

Observe that  
$$\Sigma :=\left\{x\in \overline{\Omega}: S(x)<\infty\right\} \subset \partial \Omega \,,$$
and hence the $N$-dimensional Lebesgue measure $|\overline{\Sigma}|=0.$
Indeed, if $y$ is an interior point of $\Omega$, then there exists \(r>0\) sufficiently small such that
\[
d_\Omega(x)\ge c_y>0
\quad \ \forall x \in B_{r}(y),
\]
for some constant $c_y.$
Hence, for \(u\in W_0^{1,2}(\Omega\cap B_r(y))\),

\[
\int_\Omega \frac{|u|^2}{d_\Omega^2} {\mathrm{d}}x
\le
\frac1{c_y^2}\int_\Omega |u|^2 {\mathrm{d}}x.
\]
By Poincar\'e inequality, we further have
\[
\int_\Omega |u|^2{\mathrm{d}}x
\le
Cr^2\int_\Omega |\nabla u|^2{\mathrm{d}}x
\]
for some $C>0$.
Therefore,
\[
\int_\Omega \frac{|u|^2}{d_\Omega^2} {\mathrm{d}}x
\le
\frac{C}{c_y^2}r^2\int_\Omega |\nabla u|^2{\mathrm{d}}x.
\]
This implies that $S_r(y) \geq Cc_y^{-2}\,r^{-2}.$ 
Letting \(r\to0\), we obtain
$S(y)
=
\infty,$ and hence the claim follows.

The proof of the next lemma (concentration compactness principle) follows the same lines as in \cite[Lemma 2.1]{smets} but requires a slight modification to accommodate the case $|\overline{\sum}|=0$. These required modifications were carried out in \cite[Lemma 4.40]{Hou} (see also \cite[Section 6]{DY}). 
 
\begin{lem}\label{concentration} 
Let $\Omega$ be a bounded domain in $\mathbb{R}^N$ and $(u_n)$ be a bounded sequence in $W_0^{1,2}(\Omega)$. Assume that  
$u_n \rightharpoonup u
$ weakly in $W_0^{1,2}(\Omega)$ up to a subsequence, and
$|\nabla u_n-\nabla u|^2 \rightharpoonup \mu$, 
$\frac{|u_n-u|^2}{d_{\Omega}^{2}} \rightharpoonup \nu$ in the space of positive finite measures. 
Then
\begin{enumerate}
\item[(i)] $\|\mu\| \ge S^*\|\nu\|$, where $\|\cdot\|$ represents the total variation of these measures,
    \item[(ii)] $ \overline{\lim}_{n\to\infty}
\int_{\Omega} \frac{|u_n|^2}{d_{\Omega}^2}{\mathrm{d}}x=\int_{\Omega} \frac{|u|^2}{d_{\Omega}^2}{\mathrm{d}}x+\|\nu\|,$
\item[(iii)]  $\overline{\lim}_{n\to\infty}
\int_{\Omega} |\nabla u_n|^2{\mathrm{d}}x= \int_{\Omega} |\nabla u|^2{\mathrm{d}}x+\|\mu\|.$
\end{enumerate}
\end{lem}

The following proposition compares  $S^*$ and the Hardy constant at infinity.
\begin{prop}
    \label{p-infty-estimate}
 Let $\Omega$ be a bounded domain. Then $$S^*\ge \mathbb{H}_{2,\infty}(\Omega)^2.$$
\end{prop}
\begin{proof}
Let $\lambda\ge 0$ be such that
\[
\lambda< \mathbb{H}_{2,\infty}(\Omega).
\]
By the definition of $\mathbb{H}_{2,\infty}(\Omega)$, there exists a compact set
$K\Subset\Omega$
such that 
\begin{equation}\label{Eqn1;lemm;sec3}
  \int_\Omega |\nabla u|^2\,{\mathrm{d}}x
\ge
\lambda^2
\int_\Omega \frac{|u|^2}{d_{\Omega}^2}\,{\mathrm{d}}x,\quad\forall\, u\in W_0^{1,2}(\Omega\setminus \overline{K}).  
\end{equation}
Now let \(x_0\in\partial\Omega\). Since \(K\Subset\Omega\), for sufficiently small \(r>0\),
\[
\Omega\cap B_r(x_0)\subset \Omega\setminus \overline{K}.
\]
Hence every admissible function in the definition of \(S_{r}(x_0)\) satisfies \eqref{Eqn1;lemm;sec3}. Therefore,
\[
\int_\Omega |\nabla u|^2\,{\mathrm{d}}x
\ge
\lambda^2
\int_\Omega \frac{|u|^2}{d_{\Omega}^2} \,{\mathrm{d}}x\quad\forall\, u\in W_0^{1,2}(\Omega\cap B_r(x_0)).
\]
For normalized functions (i.e.,
$\|u\|_{L^2(\Omega,\,d_\Omega^{-2})}=1$),
we obtain
\[
\int_\Omega |\nabla u|^2\,{\mathrm{d}}x\ge\lambda^2.
\]
Taking the infimum,
$S_{r}(x_0)\ge\lambda^2.$
Passing to the limit \(r\to0\), we obtain
$S(x_0)\ge\lambda^2.$
Since $x_0 \in \partial \Omega$ was arbitrary and $S(\cdot)=+\infty$ for any interior point, it follows that
\[
S^{*}
=
\inf_{x_0\in \partial \Omega} S(x_0)
\ge
\lambda^2.
\]
Finally, since the inequality holds for every \(\lambda<\mathbb{H}_{2,\infty}(\Omega)\), we conclude that
$S^*\ge\mathbb{H}_{2,\infty}(\Omega)^2.$
\end{proof}

The next proposition shows that  $\Lambda(\Omega)$ is strictly larger than the Hardy constant $\mathbb{H}_2(\Omega).$ The proof of the following result relies on Lemma \ref{concentration} and Proposition \ref{p-infty-estimate}.

\begin{prop}\label{Prop-attain} Let $\Omega$ be as in Theorem \ref{Hardy-Thm-2}. Then 
   $$\Lambda(\Omega)>\mathbb{H}_2(\Omega).$$  Moreover, $\Lambda(\Omega)$ is attained if $\mathbb{H}_{2,\infty}(\Omega)>\Lambda(\Omega).$
\end{prop}
\begin{proof}   
From the definition of $\Lambda(\Omega),$  there exists $(u_k)$ in $W^{1,2}_0(\Omega)$ with $\langle u_k,\Psi\rangle_{W^{1,2}_0(\Omega)}=0$  such that $$R_2(u_k)\to \Lambda(\Omega).$$ Since the Rayleigh quotient $R_2(\cdot)$ is homogeneous with degree $0$, by normalizing the minimizing sequence (i.e., $\|u_k\|_{L^2(\Omega,d_{\Omega}^{-2})}=1$) and calling it again $(u_k)$, we get $\|\nabla u_k \|_{L^2(\Omega)} \to \Lambda(\Omega).$ Thus, $(u_k)$ is a bounded sequence in $W^{1,2}_0(\Omega)$, and hence, it has a weakly convergent subsequence, call it again $u_k\rightharpoonup v_0$ in $W^{1,2}_0(\Omega)$. Due to weak lower semicontinuity, we have $\|\nabla v_0\|_{L^2(\Omega)} \leq \Lambda(\Omega)$. Also, observe that
 $$\int_\Omega\nabla \Psi\cdot\nabla v_0\,{\mathrm{d}}x=\lim_{k\to\infty}\int_\Omega\nabla \Psi\cdot\nabla u_k\,{\mathrm{d}}x=0,$$ i.e., $\langle v_0,\Psi\rangle_{W^{1,2}_0(\Omega)}=0.$ Applying Lemma \ref{concentration}, we get 
\begin{align}
 \|\mu\|&\ge S^*\|\nu\|\label{Eqn1:prop-3},\\
    1&=\int_{\Omega} \frac{|v_0|^2}{d_{\Omega}^2}\,{\mathrm{d}}x+\|\nu\|\label{Eqn2:prop-3},\\ \Lambda(\Omega)^2&= \int_{\Omega} |\nabla v_0|^2\,{\mathrm{d}}x+\|\mu\|\label{Eqn3:prop-3}.
\end{align} 
From the definition of $\Lambda(\Omega),$ 
\begin{align}
    \int_{\Omega} |\nabla v_0|^2\,{\mathrm{d}}x\label{Eqn4:prop-3}&\ge \Lambda(\Omega)^2 \int_\Omega \frac{|v_0|^2}{d_{\Omega}^2}\,{\mathrm{d}}x\\&\ge \mathbb{H}_2(\Omega)^2  \int_\Omega \frac{|v_0|^2}{d_{\Omega}^2}\,{\mathrm{d}}x\label{Eqn5:prop-3}.
\end{align}
Since
$S^*\ge \mathbb{H}_{2,\infty}(\Omega)^2$ (Proposition \ref{p-infty-estimate}) and spectral gap $\Gamma_2(\Omega)=\mathbb{H}_{2,\infty}(\Omega)-\mathbb{H}_2(\Omega)>0,$ 
\begin{equation}\label{eqn6:prop-3}
 S^*>\mathbb{H}_2(\Omega)^2 \,.   
\end{equation}

Now we prove
$\Lambda(\Omega)>\mathbb{H}_2(\Omega)$. Suppose  $\Lambda(\Omega)= \mathbb{H}_2(\Omega).$
Combining \eqref{Eqn1:prop-3}, \eqref{Eqn2:prop-3}, \eqref{Eqn3:prop-3}, and \eqref{Eqn5:prop-3}, we obtain $$\Lambda(\Omega)^2 \ge \mathbb{H}_2(\Omega)^2+\|\nu\|(S^*-\mathbb{H}_2(\Omega)^2).$$
Then from the above inequality, 
\begin{equation*}
    \|\nu\|(S^*-\mathbb{H}_2(\Omega)^2)\le 0.
\end{equation*} 
 and hence \eqref{eqn6:prop-3}  implies
 $\nu=0.$ Consequently, \eqref{Eqn2:prop-3} gives $$\int_\Omega \frac{|v_0|^2}{d_{\Omega}^2}\,{\mathrm{d}}x=1 \,.$$ 
Recall that (see \eqref{Eqn3:prop-3}) $\mathbb{H}_2(\Omega)^2 = \Lambda(\Omega)^2 \ge \int_\Omega|\nabla v_0|^2\,{\mathrm{d}}x.$ 
Thus, $$R_2(v_0)\le \mathbb{H}_2(\Omega).$$ Hence, $R_2(v_0)= \mathbb{H}_2(\Omega),$ i.e., $v_0$ is a minimizer for $\mathbb{H}_2(\Omega).$ Therefore, $v_0=c\,\Psi$ for some $c \neq 0$. This is a contradiction, as $\langle v_0,\Psi\rangle_{W^{1,2}_0(\Omega)}=0.$ Hence $\Lambda(\Omega)>\mathbb{H}_2(\Omega).$ 

Next, we prove that $\Lambda(\Omega)$ is attained if $\mathbb{H}_{2,\infty}(\Omega)>\Lambda(\Omega)$. Combining \eqref{Eqn1:prop-3}, \eqref{Eqn2:prop-3}, \eqref{Eqn3:prop-3}, and \eqref{Eqn4:prop-3},  we obtain
\begin{equation*}
    0\ge \|\nu\|(S^*-\Lambda(\Omega)).
\end{equation*}
Since
$S^*\ge \mathbb{H}_{2,\infty}(\Omega)$ and $\mathbb{H}_{2,\infty}(\Omega)>\Lambda(\Omega),$ the above inequality implies $\nu=0 $,
which gives $\int_{\Omega} \frac{|v_0|^2}{d_{\Omega}^2}\,{\mathrm{d}}x=1$. Also, recall that $\langle v_0,\Psi\rangle_{W^{1,2}_0(\Omega)}=0$, $\|\nabla v_0\|_{L^2(\Omega)} \leq \Lambda(\Omega)$.
This implies that $\Lambda(\Omega)$ is attained at $v_0$. 
\end{proof}
\begin{proof}[Proof of Theorem \ref{Hardy-Thm-2}]
    
We divide the three parts of the theorem as follows.

{\bf Stability.} Without loss of generality, we can assume that
$\|u\|_{L^2(\Omega,\,d_\Omega^{-2})}=1$ to establish \eqref{Hardy-2}. 
Choose such $u$ arbitrary. Since $\mathcal{M}_{\mathrm{Hardy}}$ is a closed subspace of the Hilbert space $W^{1,2}_0(\Omega)$, by the Hilbert projection theorem, $D_H(u)=\inf_{c\in\mathbb{R}} \|\nabla( u- c\,\Psi)\|_{L^2(\Omega)}$ is uniquely attained \cite[Theorem 7.1.2, Remark 7.1.8]{kesh-f}, i.e., there exists a unique $a_0=a_0(u)\in\mathbb{R}$ such that $$D_H(u) =\| \nabla(u- a_0\Psi)\|_{L^2(\Omega)} \,.$$ Moreover, we can decompose 
\begin{equation}\label{explicit}
u=a_0\Psi+D_H(u)
\,v\,,    
\end{equation}
 where  $\langle v,\Psi\rangle_{W^{1,2}_0(\Omega)}=0$ with $\|\nabla v\|_{L^2(\Omega)}=1$. 
Thus,
\begin{align}
\int_\Omega|\nabla u|^2{\mathrm{d}}x&=a_0^2\int_\Omega |\nabla \Psi|^2\,{\mathrm{d}}x+2a_0\,D_H(u)\int_\Omega\nabla \Psi\cdot\nabla v\,{\mathrm{d}}x+D_H(u)^2\int_\Omega|\nabla v|^2{\mathrm{d}}x\nonumber    \\& =a_0^2+D_H(u)^2\label{H-p=2}\,,
\end{align}
where we used $\|\nabla \Psi\|_{L^2(\Omega)}=1.$
Moreover,
\begin{align*}
\int_\Omega\frac{|u|^2}{d_\Omega^2}{\mathrm{d}}x&=a_0^2\int_\Omega\frac{|\Psi|^2}{d_\Omega^2}{\mathrm{d}}x+2a_0\,D_H(u)\int_\Omega\frac{\Psi\, v}{d_\Omega^2}{\mathrm{d}}x+D_H(u)^2\int_\Omega\frac{|v|^2}{d_\Omega^2}{\mathrm{d}}x\\&= \frac{a_0^2}{\mathbb{H}_2(\Omega)^2}+2a_0D_H(u)\int_\Omega\frac{\Psi \,v}{d_\Omega^2}{\mathrm{d}}x+D_H(u)^2\int_\Omega\frac{|v|^2}{d_\Omega^2}{\mathrm{d}}x.  
\end{align*}
Since $\langle v,\Psi\rangle_{W^{1,2}_0(\Omega)}=0$ and $\Psi$ solves \eqref{EVP}, it follows that $$\mathbb{H}_2(\Omega)^2\int_\Omega\frac{\Psi\,v}{d_\Omega^2}{\mathrm{d}}x =\int_\Omega \nabla \Psi\cdot\nabla v\,{\mathrm{d}}x=0.$$ Combining the above two equalities, we obtain $$\mathbb{H}_2(\Omega)^2\int_\Omega\frac{|u|^2}{d_\Omega^2}{\mathrm{d}}x=a_0^2+D_H(u)^2\mathbb{H}_2(\Omega)^2\int_\Omega\frac{|v|^2}{d_\Omega^2}{\mathrm{d}}x \,.$$ Consequently, using \eqref{H-p=2} we get
   $$ \delta_H(u)=\int_\Omega|\nabla u|^2{\mathrm{d}}x- \mathbb{H}_2(\Omega)^2\int_\Omega\frac{|u|^2}{d_\Omega^2}\,{\mathrm{d}}x= D_H(u)^2\left(1-\mathbb{H}_2(\Omega)^2\int_\Omega\frac{|v|^2}{d_\Omega^2}{\mathrm{d}}x\right).$$
Observe that
$$ \Lambda(\Omega)^2=\inf_{\langle u,\Psi\rangle_{W^{1,2}_0(\Omega)}=0}\frac{\int_\Omega |\nabla u|^2{\mathrm{d}}x}{\int_\Omega\frac{|u|^2}{d_\Omega^2}{\mathrm{d}}x}\le \frac{\int_\Omega |\nabla v|^2{\mathrm{d}}x}{\int_\Omega\frac{|v|^2}{d_\Omega^2}{\mathrm{d}}x}\le  \frac{1}{\int_\Omega\frac{|v|^2}{d_\Omega^2}{\mathrm{d}}x}.$$ This implies that $\int_\Omega\frac{|v|^2}{d_\Omega^2}{\mathrm{d}}x\le \frac{1}{\Lambda(\Omega)^2}.$ Therefore, $$\delta_H(u)\ge D_H(u)^2\left(1-\frac{\mathbb{H}_2(\Omega)^2}{\Lambda(\Omega)^2}\right).$$
This proves the gradient stability of the $L^2$-Hardy inequality in $\Omega.$

{\bf Sharpness.} For the second part of the theorem, we define $$\mathbb{B}:=\inf_{u\in W_0^{1,2}(\Omega)}\frac{\delta_H(u)}{D_H(u)^2}.$$ 
From the first part, we see that
\begin{equation}\label{eqBest-const}
   \mathbb{B}\ge \left(1-\frac{\mathbb{H}_2(\Omega)^2}{\Lambda(\Omega)^2}\right). 
\end{equation}
 From the definition of $\Lambda(\Omega)$, there exists $(u_k)$ with $\langle u_k,\Psi\rangle_{W^{1,2}_0(\Omega)}=0$ such that 
 \begin{equation}\label{eqn:assumption}
 \int_\Omega\frac{u_k^2}{d_\Omega^2}{\mathrm{d}}x=1,\quad\int_\Omega|\nabla u_k|^2{\mathrm{d}}x\to \Lambda(\Omega)^2.    
 \end{equation}
  Consider $\phi_k:=\Psi+u_k,\;k\in \mathbb{N}.$ 
Observe that \begin{align}
\delta_H(\phi_k)&=\frac{1}{\|\phi_k\|^2_{L^2(\Omega,\,d_\Omega^{-2})}}\Big[\int_\Omega|\nabla \Psi|^2\,{\mathrm{d}}x+\int_\Omega|\nabla u_k|^2\,{\mathrm{d}}x-\mathbb{H}_2(\Omega)^2\int_\Omega\frac{|\Psi|^2}{d_\Omega^2}{\mathrm{d}}x-\mathbb{H}_2(\Omega)^2\int_\Omega\frac{|u_k|^2}{d_\Omega^2}{\mathrm{d}}x\Big]\nonumber\\&=\frac{1}{\|\phi_k\|^2_{L^2(\Omega,\,d_\Omega^{-2})}}\Big[\int_\Omega|\nabla u_k|^2\,{\mathrm{d}}x-\mathbb{H}_2(\Omega)^2\int_\Omega\frac{u_k^2}{d_\Omega^2}{\mathrm{d}}x\Big]\label{com--1}.    
\end{align}
 On the other hand,
 \begin{align}
D_H(\phi_k)^2=&\frac{1}{\|\phi_k\|^2_{L^2(\Omega,\,d_\Omega^{-2})}}\inf_{c\in\mathbb{R}}\int_\Omega|(1-c)\nabla \Psi+\nabla u_k|^2\,{\mathrm{d}}x \nonumber \\=&\frac{1}{\|\phi_k\|^2_{L^2(\Omega,\,d_\Omega^{-2})}}\Big[\inf_{c\in\mathbb{R}}(1-c)^2\int_\Omega|\nabla \Psi|^2\,{\mathrm{d}}x+\int_\Omega|\nabla u_k|^2\,{\mathrm{d}}x\Big]=\frac{\|\nabla u_k\|^2_{L^2(\Omega)}}{\|\phi_k\|^2_{L^2(\Omega,\,d_\Omega^{-2})}}\,.
\label{com}
 \end{align}
Now using \eqref{eqn:assumption}, we get
$$\mathbb{B}\le\lim_{k\to\infty }\frac{\delta_H(\phi_k)}{D_H(\phi_k)^2}=\lim_{k\to\infty }\frac{
\|\nabla u_k\|^2_{L^2(\Omega)}-\mathbb{H}_2(\Omega)^2
}{\|\nabla u_k\|^2_{L^2(\Omega)}}=\left(1-\frac{\mathbb{H}_2(\Omega)^2}{\Lambda(\Omega)^2}\right).$$ This together with \eqref{eqBest-const} proves the sharpness of \eqref{Hardy-2}. 

\textbf{Attainability.} Assume that $\Lambda(\Omega)$ is attained at $v\in W^{1,2}_0(\Omega)\setminus\{0\}.$ Take $u$ of the form
$u=c \Psi+v$ for some arbitrary $c \in \mathbb{R}$. Now
a similar computation as before (cf. \eqref{com--1} and \eqref{com}) gives 
\begin{align*}
    \delta_H(u) &=\frac{1}{\|u\|^2_{L^2(\Omega,\,d_\Omega^{-2})}}\Big[\|\nabla v\|^2_{L^2(\Omega)}- \mathbb{H}_2(\Omega)^2\|v\|^2_{L^2(\Omega,\,d_\Omega^{-2})} \Big] \\&=\frac{1}{\|u\|^2_{L^2(\Omega,\,d_\Omega^{-2})}}\Big[\|\nabla v\|^2_{L^2(\Omega)}- \frac{\mathbb{H}_2(\Omega)^2}{\Lambda(\Omega)^2}\|\nabla v\|^2_{L^2(\Omega)} \Big]  =  
    \left( 1- \frac{\mathbb{H}_2(\Omega)^2}{\Lambda(\Omega)^2}   \right)\frac{\|\nabla v\|^2_{L^2(\Omega)}}{\|u\|^2_{L^2(\Omega,\,d_\Omega^{-2})}} \;, 
\end{align*}
 $$D_H(u)^2=\frac{\|\nabla v\|^2_{L^2(\Omega)}}{\|u\|^2_{L^2(\Omega,\,d_\Omega^{-2})}}\,,$$ which shows that the equality occurs in \eqref{Hardy-2} for $u.$
Now, conversely assume that the equality occurs in \eqref{Hardy-2} for some $u \in W^{1,2}_0(\Omega)\setminus\{0\}$. Observe that $u$ can be decomposed as $u=a_0 \Psi+v$, for some $a_0=a_0(u) \in \mathbb{R}$ and $\langle v,\Psi\rangle_{W^{1,2}_0(\Omega)}=0$ (see \eqref{explicit}). Then we have as above
$$
    \delta_H(u) =\frac{1}{\|u\|^2_{L^2(\Omega,\,d_\Omega^{-2})}}\Big[\|\nabla v\|^2_{L^2(\Omega)}- \mathbb{H}_2(\Omega)^2\|v\|^2_{L^2(\Omega,\,d_\Omega^{-2})} \Big] \quad\text{and}\quad  D_H(u)^2=\frac{\|\nabla v\|^2_{L^2(\Omega)}}{\|u\|^2_{L^2(\Omega,\,d_\Omega^{-2})}}
    \,.
$$ Since  $\delta_H(u)=\left(1-\frac{\mathbb{H}_2(\Omega)^2}{\Lambda(\Omega)^2}\right)D_H(u)^2,$
it follows that $\|\nabla v\|_{L^2(\Omega)}= \Lambda(\Omega) \|v\|_{L^2(\Omega,\,d_\Omega^{-2})}.$ That is $\Lambda(\Omega)$ is attained at $v,$ completing the proof. 
\end{proof}

\subsection{The nonlinear case: $p\neq2$}
In this subsection, we investigate the gradient stability of the $L^p$-Hardy inequality for $p\ge 2$ and $p< 2$ separately. The following reduction lemma will be used in both cases. The proof follows from the same arguments employed in the proof of Lemma \ref{Thm:Aux}, and hence we omit the details.
\begin{lem}
    \label{ThmHardy:Aux} Let $p>1$ and $\Omega$ be a bounded  $C^{1,\gamma}$-domain in $\mathbb{R}^N$ with $\Gamma_p(\Omega)>0$. Then there exist positive constants $C_1=C_1(p)$ and $C_2=C_2(p,N,\Omega)$ such
that 
\begin{equation}\label{eqn-p-case-aux}
  \left(  \frac{\|\nabla u-\nabla v\|_{L^p}}{\|u\|_{L^p(\Omega,\,d_\Omega^{-p})}}\right)^\alpha\le C_1\delta_H(u)+C_2 \frac{\| u- v\|_{L^p(\Omega,\,d_\Omega^{-p})}}{\|u\|_{L^p(\Omega,\,d_\Omega^{-p})}}.
\end{equation}
for all $u \in W^{1,p}_0(\Omega)$ and for any $v \in \mathcal{M}_{\mathrm{Hardy}}$ with $\|u\|_{L^p(\Omega,\,d_\Omega^{-p})} = \|v\|_{L^p(\Omega,\,d_\Omega^{-p})},$ 
where $\alpha=p'$ when $p<2,$ and  $\alpha=p$ when $p\ge 2.$ 
\end{lem}

\subsubsection{The case: $p\geq2$}
First, we establish the following `weak' stability result, which eventually helps us to prove the desired quantitative stability with respect to the distance $D_H(\cdot)$.
\begin{thm}\label{Hardy-Thm-p} Let $p\ge 2$ and $\Omega$ be a bounded  $C^{1,\gamma}$ domain in $\mathbb{R}^N$ such that $\Gamma_p(\Omega)>0.$ Then there exists a positive constant $C=C(p,N,\Omega)$ such that \[ \delta_H(u)\ge C \eta (u) \qquad \forall u\in W^{1,p}_0(\Omega),\] where
\begin{equation}\label{Eqn-Hardy-p}
    \eta(u):= \inf\left\{\frac{\|u-v\|^{p}_{L^p(\Omega,\,d_\Omega^{-p})}}{\|u\|^{p}_{L^p(\Omega,\,d_\Omega^{-p})}}:v\in \mathcal{M}_{\mathrm{Hardy}},\,\|v\|_{L^p(\Omega,\,d_\Omega^{-p})}=\|u\|_{L^p(\Omega,\,d_\Omega^{-p})}\right\}\,.
\end{equation}
\end{thm}

We recall two lemmas, which will be used to prove Theorem \ref{Hardy-Thm-p}. The first one is the following  weighted Poincar\'e inequality  \cite[Theorem 5.3]{Poincare2025}. See also \cite[Theorem 1.12]{Poincare_2011} for a proof.
\begin{lem}[Weighted Poincar\'e inequality]\label{poincare}
Let $\Omega$ be a bounded $C^{1,\gamma}$-domain in $\mathbb{R}^N$, $p>1,$ and  $\beta>- \frac{1}{p}.$ Then there exists a positive constant $C=C(\Omega,\beta,p,N)$ such that 
\begin{equation*}
    \int_\Omega|\nabla u|^p \,d_\Omega^{(\beta+1)p}\,{\mathrm{d}}x\ge C\int_\Omega|u(x)-u_{\Omega,\beta}|^p \,d_\Omega^{\beta p}\, {\mathrm{d}}x\qquad \forall \,u\in W^{1,p}(\Omega)\cap C_c(\Omega) \,,
\end{equation*}
where $u_{\Omega,\beta}$ is the weighted average of $u$ over $\Omega$, i.e.,
$u_{\Omega,\beta} := \left(\int_\Omega u\, d_\Omega^\beta\, {\mathrm{d}}x\right)\left(\int_\Omega d_\Omega^\beta\, {\mathrm{d}}x\right)^{-1}.$
\end{lem}

The second lemma concerns the equivalence of the {\it energy functional} and the {\it simplified energy functional}. We will replace the role of the following energy functional associated to \eqref{EVP}
$$Q_{p}(u)=\int_\Omega[|\nabla u|^p- \mathbb{H}_p(\Omega)^p\frac{|u|^p}{d_{\Omega}^p}] \ {\mathrm{d}}x\qquad\,\forall\, u\in W^{1,p}(\Omega) \cap C_c(\Omega)$$
by the simplified energy functional given by
$$Q_{p}^{sim}(u)=\int_\Omega w^2|\nabla u|^2\left(u|\nabla w|+w|\nabla u|\right)^{p-2}{\mathrm{d}}x \quad\,\forall\, u\in W^{1,p}(\Omega) \cap C_c(\Omega) \,,$$
where $w  \in W^{1,p}_\mathrm{loc}(\Omega)\cap C(\Omega)$ is a positive weak solution of the following  equation:
\begin{equation*}
    \mathcal{Q}'_{p}[\phi]:=-\Delta_p \phi- \mathbb{H}_p(\Omega)^p \frac{|\phi|^{p-2}\phi}{d_{\Omega}^p}=0 \ \,\text{in} \, \ \Omega \,.
\end{equation*}

\begin{lem}{\mbox{(\cite[Lemma 2.2]{Yehuda2008} and \cite[Lemma 3.4]{Yehuda2015})}} \label{prop-1}
    Let $p\in (1,\infty)$ and $w\in W^{1,p}_\mathrm{loc}(\Omega)\cap C(\Omega)$ be a positive weak solution of $\mathcal{Q}_{p}^\prime[\phi]=0$. Then $\exists C_1,C_2>0$ (depending only on $p$) such that  $$ C_1 Q_{p}( u) \leq Q_{p}^{sim}(\frac{u}{w}) \leq C_2 Q_{p}( u) \qquad (\mbox{in short,} \ Q_{p}( u) \asymp Q_{p}^{sim}(\frac{u}{w})) $$
    for all non-negative $u\in W^{1,p}(\Omega) \cap C_c(\Omega)$.
\end{lem}

We will use the above lemmas to prove  Theorem \ref{Hardy-Thm-p}. Additionally, the boundary behavior of the extremal $\Psi$ will play a key role in our proof.
 \begin{proof}[Proof of Theorem \ref{Hardy-Thm-p}]
Due to the homogeneity-free formulation of \eqref{Eqn-Hardy-p}, it is enough to prove 
\begin{align} \label{Eq:main_req}
    \delta_H(u)\gtrsim \inf\left\{\|u-v\|^{p}_{L^p(\Omega,\,d_\Omega^{-p})}:v\in \mathcal{M}_{\mathrm{Hardy}},\,\|v\|_{L^p(\Omega,\,d_\Omega^{-p})}=\|u\|_{L^p(\Omega,\,d_\Omega^{-p})}\right\}
\end{align}
for all $u\in W^{1,p}_0(\Omega)$ with $\|u\|_{L^p(\Omega,\,d_\Omega^{-p})}=1.$ Further, since $W^{1,p}(\Omega) \cap C_c(\Omega)$ is dense in $W^{1,p}_0(\Omega),$ it is sufficient to prove \eqref{Eq:main_req}
for all $u\in W^{1,p}(\Omega) \cap C_c(\Omega)$ with $\|u\|_{L^p(\Omega,\,d_\Omega^{-p})}=1.$ Choose such $u$ arbitrarily.

\noindent We divide our proof into two steps. 

{\bf Step-1.} We will show that 
\begin{align} \label{Eq:req1}
    \delta_H(u)\gtrsim\int_\Omega \Psi^p\left|\nabla \left(\frac{u}{\Psi}\right)\right|^p{\mathrm{d}}x.
\end{align}
 Applying Lemma \ref{prop-1} and using the fact that $p\ge 2$, we obtain 
\begin{align}
  \delta_H(u)=  \delta_H(|u|)=Q_{p}(|u|) & \ge C \int_\Omega \Psi^2\left|\nabla \left(\frac{|u|}{\Psi}\right) \right|^2\left(\frac{|u|}{\Psi}|\nabla \Psi|+\Psi\left|\nabla \left(\frac{|u|}{\Psi}\right) \right|\right)^{p-2}{\mathrm{d}}x \nonumber \\
  &\ge C \int_\Omega \Psi^p\left|\nabla \left(\frac{|u|}{\Psi}\right)\right|^p{\mathrm{d}}x \label{p2} \,
\end{align}
for some $C=C(p)>0,$ concludes \eqref{Eq:req1}.

{\bf Step-2.} We will show that 
\begin{equation}\label{Eqn1-thm-p-Hardy}
    \int_\Omega \Psi^p\left|\nabla \left(\frac{u}{\Psi}\right)\right|^p{\mathrm{d}}x\gtrsim\inf\left\{\|u-v\|^{p}_{L^p(\Omega,\,d_\Omega^{-p})}:v\in \mathcal{M}_{\mathrm{Hardy}}\right\}.
\end{equation}
It is known that \cite[Section 4,5]{Pinch19} (see also \cite[Theorem 1.1]{Itai} for $C^2$ domain) there exists $C_1=C_1(p,N,\Omega)>0$ such that
$$C_1^{-1}d_\Omega^\alpha(x)\le \Psi(x)\le C_1d_\Omega^\alpha(x)\ \ \text{in} \ \Omega \ \ \mbox{for some} \  \ \alpha\in \left(1-\frac{1}{p},1\right) \,,$$
where $\alpha$ satisfies $ (p-1)\alpha^{p-1}(1-\alpha)=\mathbb{H}_p(\Omega)^p$. Thus, apply Lemma \ref{poincare} with $\beta=\alpha-1>-\frac{1}{p}$ (thanks to the range of $\alpha$) to get
\begin{align*}
     \int_\Omega \Psi^p\left|\nabla u\right|^p{\mathrm{d}}x&\ge C_1^{-p} \int_\Omega|\nabla u|^p d_\Omega^{\alpha p}\,{\mathrm{d}}x\\& \ge C_1^{-p}C_2 \int_\Omega|u(x)-u_{\Omega,\alpha}|^p d_\Omega^{(\alpha-1) p}\, {\mathrm{d}}x \ge C_1^{-p}C_2 \inf_{c\in \mathbb{R}}\int_\Omega|u(x)-c|^p d_\Omega^{(\alpha-1) p}\, {\mathrm{d}}x \,,
\end{align*}
for some constant $C_2=C_2(p,N,\Omega)>0.$
 Now, replace $u$ with $\frac{u}{\Psi}$ to obtain 
\begin{align*}
     \int_\Omega \Psi^p\left|\nabla \left(\frac{u}{\Psi}\right)\right|^p{\mathrm{d}}x&\ge C_1^{-p}C_2\inf_{c\in \mathbb{R}}\int_\Omega\left|\frac{u}{\Psi}-c\right|^pd_\Omega^{(\alpha-1) p}\,{\mathrm{d}}x\\&\ge C_1^{-2p}C_2 \inf_{c\in \mathbb{R}}\int_\Omega\left|\frac{u}{\Psi}-c\right|^p\frac{\Psi^p}{d_\Omega^{p}}
    \,{\mathrm{d}}x =C_1^{-2p}C_2\inf_{c\in \mathbb{R}}\int_\Omega\frac{\left|u-c\Psi\right|^p}{d_\Omega^p}{\mathrm{d}}x
\end{align*}
i.e., \eqref{Eqn1-thm-p-Hardy} holds. Now, \eqref{Eq:main_req} follows by combining \eqref{Eq:req1} and \eqref{Eqn1-thm-p-Hardy}, and applying Theorem
\ref{Thm:A}. 
This completes our proof.
\end{proof}

Now we will prove the desired stability of the $L^p$-Hardy inequality for $p \geq 2$ with respect to the distance $D_H(\cdot)$. The proof follows by using the reduction lemma (Lemma \ref{ThmHardy:Aux}) together with Theorem \ref{Hardy-Thm-p}.
\begin{thm}\label{Hardy-Thm-grad}  Let $p\ge 2$ and $\Omega$ be a bounded  $C^{1,\gamma}$-domain in $\mathbb{R}^N$ such that $\Gamma_p(\Omega)>0.$ Then there exists a positive constant $C=C(p,N,\Omega)$  
such that
\begin{equation}\label{Hardy-2=grad}
  \delta_H(u)\ge C \min\left\{D_H(u)^{p},\,D_H(u)^{p^2}\right\}, \quad\forall\,u\in W^{1,p}_0(\Omega).  
\end{equation}
\end{thm}
 \begin{proof}
By Theorem \ref{Hardy-Thm-p} we get a positive constant  $C=C(p,N,\Omega)$ such that
\begin{equation}
    \label{eqn-Hardy-grad}\delta_H(u)\ge C\,\eta(u)^p\,\quad \forall\,u\in W_0^{1,p}(\Omega)\,.
\end{equation}
If $\eta(u) =0$, then it follows that $u \in \mathcal{M}_{\mathrm{Hardy}}$. Hence,  $D_H(u)=0$. This implies that Theorem \ref{Hardy-Thm-grad} holds trivially. Now assume that $\eta(u) >0$. Let $v\in \mathcal{M}_{\mathrm{Hardy}}$ be such that $$\|u\|_{L^p(\Omega,\,d_\Omega^{-p})} = \|v\|_{L^p(\Omega,\,d_\Omega^{-p})}\quad\text{and}\quad\frac{\| u- v\|_{L^p(\Omega,\,d_\Omega^{-p})}}{\|u\|_{L^p(\Omega,\,d_\Omega^{-p})}}\le 2\eta(u).$$  Then apply  Lemma \ref{ThmHardy:Aux} (for $p\ge 2$) and \eqref{eqn-Hardy-grad} to obtain
 \begin{align*}
    D_H(u)^p\le\left(\frac{\|\nabla u-\nabla v\|_{L^p(\Omega)}}{\|u\|_{L^p(\Omega,\,d_\Omega^{-p})}}\right)^p&\le C_1\delta_H(u)+C_2\frac{\| u- v\|_{L^p(\Omega,\,d_\Omega^{-p})}}{\|u\|_{L^p(\Omega,\,d_\Omega^{-p})}}\\&\le C_1\delta_H(u)+2C_2\eta(u) \le C_1\delta_H(u)+2C_2C^{-1/p}\delta_H(u)^\frac{1}{p}\\& \qquad \qquad \qquad \qquad  \qquad \le C_3\max\{\delta_H(u),\,\delta_H(u)^\frac{1}{p}\}, 
 \end{align*}
for some constant $C_3=C_3(p,N,\Omega)>1$.  Hence \eqref{Hardy-2=grad} follows from the above inequality and the fact that $a^p\le C_3\max\{b,b^{1/p}\} , \;a,b\ge 0$ implies $b\ge C_3^{-p}\min\{a^p,a^{p^2}\}$.
  \end{proof}

\subsubsection{The case: $1<p<2$} In this case, the previous ideas fail since we can not find the lower bound of the simplified energy as in \eqref{p2}. Nevertheless, we manage to prove a stability result which states the following:

\begin{thm}\label{Hardy-Thm-grad<2}  Let $p\in (1,2)$ and $\Omega$ be a bounded  $C^{1,\gamma}$-domain in $\mathbb{R}^N$ with $\Gamma_p(\Omega)>0.$ Then there exists a positive constant $C=C(p,N,\Omega)$ 
such that for every $u\in W^{1,p}_0(\Omega),$
\begin{equation}\label{Hardy-2=grad=2}
    \delta_H(u)\ge C \min\left\{D_H(u)^{p'},\,D_H(u)^{2p'}\right\} \,.  
\end{equation}

\end{thm}

To prove Theorem \ref{Hardy-Thm-grad<2}, we will use the following lemma, which is actually a weaker version of the theorem and an analogue of Theorem \ref{Hardy-Thm-p}. 

\begin{lem}\label{Thm<2}
     Let $p\in (1,2)$ and $\Omega$ be a bounded  $C^{1,\gamma}$ domain in $\mathbb{R}^N$ with $\Gamma_p(\Omega)>0.$ Then there exists a positive constant $C=C(p,N,\Omega)$ such that for every $u\in W^{1,p}_0(\Omega),$
\begin{equation}\label{Eqn-Hardy-p=2}
  \frac{\|\nabla u\|_{L^p(\Omega)}^p-\mathbb{H}_p(\Omega)^p\|u\|^p_{L^p(\Omega,\,d_\Omega^{-p})}}{\|u\|^p_{L^p(\Omega,\,d_\Omega^{-p})}}\ge C\, \min\left\{\eta(u)^p,\eta(u)^2\right\} \,,
\end{equation}
where $\eta(\cdot)$ is the distance function given by \eqref{Eqn-Hardy-p}.
\end{lem}

\begin{proof}
It can be seen that $\eta(tu)=\eta(u)$ for all $t\ne 0.$ Moreover, $W^{1,p}(\Omega) \cap C_c(\Omega)$ is dense in $W^{1,p}_0(\Omega).$ Thus, it is enough to prove \eqref{Eqn-Hardy-p=2} for all $u\in W^{1,p}(\Omega) \cap C_c(\Omega)$ satisfying $\|u\|_{L^p(\Omega,\,d_\Omega^{-p})}=1.$ Choose such $u$ arbitrarily.
Note that in this case, the left hand side of  \eqref{Eqn-Hardy-p=2} is
$$Q_{p}(u)=\int_\Omega[|\nabla u|^p- \mathbb{H}_p(\Omega)^p\frac{|u|^p}{d_{\Omega}^p}] \ {\mathrm{d}}x\,.
$$ 
Apply Lemma \ref{prop-1} to obtain 
\begin{align*}
  Q_{p}(u)= Q_{p}(|u|)  \ge C \int_\Omega \Psi^2\left|\nabla \left(\frac{|u|}{\Psi}\right) \right|^2\left(\frac{|u|}{\Psi}|\nabla \Psi|+\Psi\left|\nabla \left(\frac{|u|}{\Psi}\right) \right|\right)^{p-2}{\mathrm{d}}x,
\end{align*}
for some $C=C(p)>0.$ Set $\tilde{u}=\frac{|u|}{\Psi}.$ One can verify that for $p\in(1,2)$ $$(|a|+|b|)^{p-2}\ge \min\{|a|^{p-2},|b|^{p-2}\}\quad \forall \,a,b\in \mathbb{R}^N.$$
Thus, use the above inequality for $a=\tilde{u}\nabla \Psi$ and $b=\Psi\nabla \tilde{u}$ to get
\begin{align}
     Q_{p}(u) &\ge C \int_\Omega\Psi^2|\nabla \tilde{u}|^2\,\min\left\{ |\tilde{u}|^{p-2} |\nabla \Psi|^{p-2},\;\Psi^{p-2}|\nabla \tilde{u}|^{p-2}\right\}{\mathrm{d}}x\nonumber\\&= \int_{A^c\cap \Omega} \Psi^p|\nabla \tilde{u}|^p {\mathrm{d}}x+ \int_{A}\Psi^2|\nabla \tilde{u}|^2|\nabla \Psi|^{p-2}|\tilde{u}|^{p-2}{\mathrm{d}}x\,,\label{eq1:p=2}
\end{align}
where $$A:=\{x\in \Omega: \Psi^{p-2}|\nabla \tilde{u}|^{p-2}>|\nabla \Psi|^{p-2}|\tilde{u}|^{p-2}\}.$$ Now, use H\"older inequality for conjugate pair $\left(\frac{2}{p},\frac{2}{2-p}\right)$ with $f=(\Psi|\nabla \tilde{u}|)^p(|\tilde{u}||\nabla \Psi|)^\frac{p(p-2)}{2}$ and $g=(|\nabla \Psi||\tilde{u}|)^\frac{p(2-p)}{2}$ to obtain 
\begin{equation}\label{Eq2:p=2}
    \int_A\Psi^p|\nabla \tilde{u}|^p{\mathrm{d}}x \le \left(\int_{A}\Psi^2|\nabla \tilde{u}|^2|\nabla \Psi|^{p-2}|\tilde{u}|^{p-2}{\mathrm{d}}x\right)^\frac{p}{2}\left(\int_{A}|\nabla \Psi|^p|\tilde{u}|^{p}{\mathrm{d}}x\right)^\frac{2-p}{2}.
\end{equation}
Moreover, using Lemma \ref{regularity} we have the  estimate:  \begin{equation}\label{regularity-lemma}
    \int_{A}|\nabla \Psi|^p|\tilde{u}|^{p}{\mathrm{d}}x\le C_1  \int_{A}\frac{\Psi^{ p}|\tilde{u}|^{p}}{ d_\Omega^{p}}{\mathrm{d}}x=C_1\int_{A}\frac{|u|^{p}}{ d_\Omega^{p}}{\mathrm{d}}x= C_1
\end{equation}
for some constant $C_1=C_1(p,N,\Omega)>0$.  Combining \eqref{eq1:p=2},  \eqref{Eq2:p=2}, and \eqref{regularity-lemma}  we get a positive constant $C_2=C_2(p,N,\Omega)$ such that
\begin{align} 
Q_p(u)&\ge C_2 \int_{A^c\cap \Omega}\Psi^p|\nabla \tilde{u}|^p{\mathrm{d}}x+ \left(\int_A\Psi^p|\nabla \tilde{u}|^p{\mathrm{d}}x\right)^\frac{2}{p} :=C_2[I_1(u)+I_2(u)^\frac{2}{p}]\label{Eq6<p},
\end{align}
Using the elementary inequality $a+b^\frac{2}{p}\ge C_3\min\{(a+b),(a+b)^\frac{2}{p}\}\; \  \forall a,b\ge0,$
 $C_3=C_3(p)>0$, we get
\begin{equation}\label{Eq7<2}
    I_1(u)+I_2(u)^\frac{2}{p}\ge C_3\min\left\{\int_{\Omega}\Psi^p|\nabla \tilde{u}|^p{\mathrm{d}}x,\left( \int_{\Omega}\Psi^p|\nabla \tilde{u}|^p{\mathrm{d}}x\right)^\frac{2}{p}\right\}.
\end{equation} 
Now, following the same arguments as in Step-2 of the proof of Theorem \ref{Hardy-Thm-p}, we can estimate
 \begin{align}
      \int_{\Omega}\Psi^p|\nabla \tilde{u}|^p{\mathrm{d}}x= \int_\Omega \Psi^p\left|\nabla \left(\frac{u}{\Psi}\right)\right|^p{\mathrm{d}}x&\ge C_4\inf\left\{\|u-v\|^{p}_{L^p(\Omega,\,d_\Omega^{-p})}:v\in \mathcal{M}_{\mathrm{Hardy}}\right\}\nonumber\\  &\ge\frac{C_4}{2} \eta(u)^p\label{Eq8<2},
 \end{align}
for some constant $C_4=C_4(p,N,\Omega)>0,$  
where we apply Theorem \ref{Thm:A} to the last inequality.
Hence \eqref{Eqn-Hardy-p=2} follows by combining \eqref{Eq6<p}, \eqref{Eq7<2}, and \eqref{Eq8<2}. 
\end{proof}

Using the above lemma, we are now ready to give a proof of Theorem \ref{Hardy-Thm-grad<2}.

\begin{proof}[Proof of Theorem \ref{Hardy-Thm-grad<2}] 
First, we verify that 
\begin{equation}\label{Eq1:m<2}
  \frac{\|\nabla u\|_{L^p(\Omega)}^p-\mathbb{H}_p(\Omega)^p\|u\|^p_{L^p(\Omega,\,d_\Omega^{-p})}}{\|u\|^p_{L^p(\Omega,\,d_\Omega^{-p})}} \le \mathbb{H}_p(\Omega)^{p-p'}\,\delta_H(u).
\end{equation}
To see this, note that for any $p\in ( 1,2),$ the function $t\mapsto t^{p'}-t^p\quad\text{for} \, t\ge 1$ is an increasing function. So, if $a\ge b\ge1,$ we have 
\begin{equation}\label{Eq2;m<2}
 a^p-b^p \leq   a^{p'}-b^{p'}  \,.
\end{equation}
Taking $a=\frac{\|\nabla u\|_{L^p(\Omega)}}{\mathbb{H}_p(\Omega)\|u\|_{L^p(\Omega,\,d_\Omega^{-p})}}$ and $b=1$ in \eqref{Eq2;m<2}, we obtain \eqref{Eq1:m<2}.
Combining \eqref{Eqn-Hardy-p=2} and \eqref{Eq1:m<2} we obtain  $\delta_H(u)\ge C\min\left\{\eta(u)^p,\eta(u)^2\right\}$ for some $C=C(p,N,\Omega)>0,$ which implies that
\begin{equation}\label{eq3<2}
     \eta(u)\le C_1\max\{\delta_H(u)^\frac{1}{2},\,\delta_H(u)^\frac{1}{p}\}\,,
\end{equation}
for some $C_1=C_1(p,N,\Omega)>0$. If $\eta(u)=0,$ then $u\in\mathcal{M}_{\mathrm{Hardy}}$. Consequently,  $D_H(u)=0.$ Hence the theorem follows trivially. Assume that $\eta(u)\neq 0$. Let $v\in \mathcal{M}_{\mathrm{Hardy}}$ be such that $$\|u\|_{L^p(\Omega,\,d_\Omega^{-p})} = \|v\|_{L^p(\Omega,\,d_\Omega^{-p})}\quad\text{and}\quad\frac{\| u- v\|_{L^p(\Omega,\,d_\Omega^{-p})}}{\|u\|_{L^p(\Omega,\,d_\Omega^{-p})}}\le 2\eta(u).$$  Then apply  Lemma \ref{ThmHardy:Aux} (for $p<2$)  and \eqref{eq3<2} to obtain
 \begin{align*}
     D_H(u)^{p'} & \le\left(\frac{\|\nabla u-\nabla v\|_{L^p(\Omega)}}{\| u\|_{L^p(\Omega,\,d_\Omega^{-p})}}\right)^{p'}\le C_2\delta_H(u)+C_2\frac{\| u- v\|_{L^p(\Omega,\,d_\Omega^{-p})}}{\|u\|_{L^p(\Omega,\,d_\Omega^{-p})}}\nonumber\\&\le C_2\delta_H(u)+2C_2\eta(u) \ \le C_2\delta_H(u)+2C_1C_2\max\{\delta_H(u)^\frac{1}{2},\,\delta_H(u)^\frac{1}{p}\}\\& \qquad \qquad \qquad \quad \qquad \quad \le C_3\max\{\delta_H(u),\delta_H(u)^\frac{1}{2}\}, 
 \end{align*}
for some $C_3=C_3(p,N,\Omega)>0$. Hence  \eqref{Hardy-2=grad=2} follows from the above inequality using the same argument as in the proof of Theorem \ref{Hardy-Thm-grad}. This completes the proof.
  \end{proof}

Instead of $D_H(\cdot)$, let us  consider the following scale invariant distance function as in \eqref{fig}
\begin{equation*}
  \tilde{D}_H(u):=\inf\left\{\frac{\|\nabla u-\nabla v\|_{L^p(\Omega)}}{\| \nabla u\|_{L^p(\Omega)}}:v\in \mathcal{M}_{\mathrm{Hardy}}\right\} \,.  
\end{equation*}
 The $L^p$-Hardy inequality infers that $D_H(u)\ge \mathbb{H}_p(\Omega) \tilde{D}_H(u)\,.$
Hence, as a corollary of Theorems \ref{Hardy-Thm-grad} and \ref{Hardy-Thm-grad<2}, we obtain the following variant of the stability inequalities \eqref{Hardy-2=grad} and \eqref{Hardy-2=grad=2}.
\begin{cor}\label{Rem-p>2} Let $p\in (1,\infty)$ and $\Omega$ be a bounded  $C^{1,\gamma}$-domain in $\mathbb{R}^N$ such that $\Gamma_p(\Omega)>0.$ Then there exists a positive constant $C=C(p,N,\Omega)$ such that 
\[
\delta_H(u)\ge C
\begin{cases}
\tilde{D}_H(u)^{2p'}
& \text{if } p<2 \\[2ex]
\tilde{D}_H(u)^{p^2}
& \text{if } p\ge 2.
\end{cases} \quad\forall\,u\in W^{1,p}_0(\Omega)
\]
\end{cor}
 
\section*{Appendix}\label{sec_appendix}

\subsection{Gradient pseudo-stability of the Morrey inequality in $\R^N$} As in the one-dimensional case (Theorem \ref{Morr-th-1d}), here we derive  a gradient pseudo-stability inequality for the Morrey inequality in $\R^N$ \eqref{Rn_morrey}. Although the result is new, we include it in the appendix, as the proof follows from a careful inspection of Cianchi's proof \cite{cian2008}, with standard modifications based on the reduction Lemma \ref{Thm:Aux:Rn}.

For weakly differentiable function $u$ in $\R^N$ with $|\operatorname{supp}u|<\infty$ and $|\nabla u|\in L^p(\R^N)$, first we  consider the following distance of $u^\star$ from the space of Morrey extremals $\mathcal{M}_{\mathrm{Morrey}}(\R^N):=\{v_{a,b,x_0}:  a\in \R,\,b\ge 0,\, x_0\in \R^N\}$:
 $$\tau(u^\star)=\inf_{a,b,x_0}\left\{\frac{\|u^\star-v_{a,b,x_0}\|_{L^\infty(\R^N)}}{\|u^\star\|_{L^\infty(\R^N)}}:  \,\|v_{a,b,x_0}\|_{L^\infty(\R^N)}=\|u^\star\|_{L^\infty(\R^N)} ,\; \operatorname{supp}v_{a,b,x_0}=\operatorname{supp}u^\star\right\}.$$
 \begin{thm}\label{Rn_sup}
     Let $p>N> 1$. Then there exists $C=C(p,N)>0$ such that 
     \begin{equation}\label{Rn_morry-sup}
         \delta_{M,\mathbb{R}^N}(u):=  \frac{\left(\kappa
|\operatorname{supp}u|^{\frac1N-\frac1p}\|\nabla u\|_{L^p(\mathbb R^N)}\right)^p
-\|u\|^p_{L^\infty(\mathbb R^N)}}{
\|u\|^p_{L^\infty(\mathbb R^N)}
}
\ge C\tau(u^\star)^p
     \end{equation}
      for every weakly differentiable function $u$ in $\R^N$ with $|\operatorname{supp}u|<\infty$ and $|\nabla u|\in L^p(\R^N).$ 
 \end{thm}

\begin{proof}
It is enough to prove \eqref{Rn_morry-sup} for every  
$u$ with
 $$\|u\|_{L^\infty(\mathbb R^N)}=1\qquad \mbox{and}\qquad|\operatorname{supp}u|=1.$$
 Given such a $u$, our aim is to find a Morrey extremal $v_{a,b,0}$ such that $\mathrm{supp} \, v_{a,b,0}= \mathrm{supp} \, u^\star = B_b(0)$ with $b=\omega_N^{-\frac{1}{N}},$ and $\|v_{a,b,0}\|_{L^{\infty}(B_b(0))} = \|u^\star\|_{L^{\infty}(B_b(0))} $, which is equal to $1$.
 
 Towards this, we follow the proof of Theorem 1.1 in \cite{cian2008}, where the quantitative stability of \eqref{Rn_morrey} is established (see \eqref{Morrey_Rn}). The author first proved (see (2.21) of \cite[Theorem 1.1]{cian2008}) that
\begin{equation}\label{Rn:E1}
    \|u^*-\phi\|_{L^\infty((0,1))} \le C \epsilon^\frac{1}{p} \left(\int_0^1\left(N\omega_N^{1/N}s^{1/N'}(-u^{*'}(s))\right)^p ds \right)^{\frac{1}{pp'}}=C \epsilon^\frac{1}{p}\|\nabla u^\star\|_{L^p(\R^N)}^\frac{1}{p'} 
\end{equation}
where $C=C(p,N)$ is a positive constant, \begin{equation}\label{4a}
           \epsilon = \kappa^p\|\nabla u\|_{L^p(\R^N)}^p-1 \ge \kappa^{p} \|\nabla u^\star\|_{L^p(\R^N)}^{p}-1 \geq \kappa \|\nabla u^\star\|_{L^p(\R^N)}-1\geq 0 \,,
     \end{equation}
   and $\phi:[0,1]\to [0,\infty)$ is defined as
 $$\phi(s)=\frac{1-s^{1-\frac{p'}{N'}}}{\Lambda^{p'}(1-\frac{p'}{N'})}\,,\qquad\, \Lambda=N\omega_N^{\frac{1}{N}} \kappa^{\frac{1}{p}}\left(\int_0^1\left(N\omega_N^{1/N}s^{1/N'}(-u^{*'}(s))\right)^p  {\mathrm{d}}s\right)^{-\frac{1}{pp'}} \,.$$   This particular choice of $\phi$ enable us to relate it to the Morrey extremals $v_{a,b,0}$ as follows $$v_{a_{\Lambda},b,0}(x) =\phi(\omega_N|x|^N)\quad\text{if}\ \omega_N|x|^N\le 1,\;\; \ \text{where} \ a_{\Lambda}=\frac{\omega_N^{1-\frac{p'}{N'}}}{\Lambda^{p'}(1-\frac{p'}{N'})}\ \text{and}\  b=\omega_N^{-\frac{1}{N}}
 \,.$$ 
It is clear that $\|u^\star-v_{a_{\Lambda},b,0}\|_{L^\infty(\R^N)} = \|u^*-\phi\|_{L^\infty((0,1))}$ (as $u^*,\phi$ is supported in $[0,1]$). Thus, using \eqref{Rn:E1} and  \eqref{4a}, we get  \begin{equation}\label{RnE2}
    \|u^\star-v_{a_{\Lambda},b,0}\|_{L^\infty(\R^N)}  \le C_1\,   \epsilon^\frac{1}{p}\left(1+\epsilon\right)^\frac{1}{p'}\,,
\end{equation}
for some $C_1=C_1(p,N)>0$.
Clearly,  $\operatorname{supp}v_{a_{\Lambda},b,0}=B_b(0)=\operatorname{supp}u^\star$, which satisfies one of our requirements. Moreover, it can be verified that 
\begin{equation}\label{Rn:E3}
    \|v_{a_{\Lambda},b,0}\|_{L^\infty(\R^N)}=\Lambda^{-p'}(1-p'/N')^{-1}=\kappa \|\nabla u^\star\|_{L^p(\R^N)} \,.
\end{equation}
 Note that $\|v_{a_{\Lambda},b,0}\|_{L^\infty(\R^N)}$ may not be equal to $\| u\|_{L^\infty(\R^N)}$, which is $1$.
Thus to satisfy our second requirement, we simply normalize $v_{a_{\Lambda},b,0}$ as $$\widetilde{v}_{a_{\Lambda},b,0}(x)=\frac{v_{a_{\Lambda},b,0}(x)}{\|v_{a_{\Lambda},b,0}\|_{L^\infty(\R^N)}},\quad x\in \R^N$$ and estimate $\|v_{a_{\Lambda},b,0}-\widetilde{v}_{a_{\Lambda},b,0}\|_{L^\infty(\R^N)}.$
Now, using \eqref{Rn:E3} and \eqref{4a} we get
\begin{equation}\label{eqn10:sup}
\|v_{a_\Lambda,b,0}-\widetilde{v}_{a_{\Lambda},b,0}\|_{L^\infty(\R^N)}=\|v_{a_{\Lambda},b,0}\|_{L^\infty(\R^N)}-1=\kappa \|\nabla u^\star\|_{L^p(\R^N)}-1\le \epsilon.  
\end{equation}
Thus, by the triangular inequality, \eqref{RnE2} and \eqref{eqn10:sup}, 
\begin{equation}\label{RN:E4}
 \|u^\star-\widetilde{v}_{a_{\Lambda},b,0}\|_{L^\infty(\R^N)}\le \|u^\star-v_{a_{\Lambda},b,0}\|_{L^\infty(\R^N)}+\|v_{a_{\Lambda},b,0}-\widetilde{v}_{a_{\Lambda},b,0}\|_{L^\infty(\R^N)}\le  C_1\epsilon^\frac{1}{p}(1+\epsilon)^\frac{1}{p'}+ \epsilon.
\end{equation}

Case 1. $\epsilon\ge 1.$ Then \eqref{4a} gives $\delta_{M,\mathbb{R}^N}(u)=\epsilon\ge 1$. Consequently, $\tau(u^\star)\le 2\le 2\delta_{M,\mathbb{R}^N}(u)$. 

Case 2. $\epsilon< 1.$ Combining \eqref{RN:E4} and \eqref{4a} we obtain $$\tau(u^\star)\le \|u^\star-\widetilde{v}_{a_{\Lambda},b,0}\|_{L^\infty(\R^N)}\le (2^\frac{1}{p'}C_1+1)\, \epsilon^\frac{1}{p}= (2^\frac{1}{p'}C_1+1) \delta_{M,\mathbb{R}^N}(u)^\frac{1}{p}.$$ Hence \eqref{Rn_morry-sup} follows by combining Case 1 and Case 2.
\end{proof}
We will be using the following reduction lemma.
\begin{lem}\label{Thm:Aux:Rn}
Let $p>N>1$. Then there exist positive constants $C_1=C_1( p)$ and $C_2=C_2(N, p)$ such that 
\begin{equation*}
  \left(  \frac{\|\nabla u^\star-\nabla v\|_{L^p(\R^N)}}{|\operatorname{supp}u^\star|^{\frac{1}{p}-\frac{1}{N}}\|u^\star\|_{L^\infty(\R^N)}}\right)^p\le C_1\delta_{M,\mathbb{R}^N}(u)+C_2 \frac{\| u^\star- v\|_{L^\infty(\R^N)}}{\|u^\star\|_{L^\infty(\R^N)}}
\end{equation*}
 for every weakly differentiable function $u$ in $\R^N$ with $|\operatorname{supp}u|<\infty$ and $|\nabla u|\in L^p(\R^N),$  and for every $v \in \mathcal{M}_{\mathrm{Morrey}}(\R^N)$ with $\|v\|_{L^\infty(\R^N)}=\|u\|_{L^\infty(\R^N)} $ and $\operatorname{supp}u^\star=\operatorname{supp}v.$ 
 
\end{lem}
\begin{proof}
   The proof follows from the same arguments employed in the proof of Lemma \ref{Thm:Aux}, and appropriate application of P\'olya--Szeg\"o inequality \eqref{polya} and the fact that  $\|u^\star\|_{L^{\infty}(\mathbb{R}^N)}=\|u\|_{L^{\infty}(\mathbb{R}^N)}$, and hence we omit the details. Also the assumption $\operatorname{supp}u^\star=\operatorname{supp}v$ will be crucial in obtaining the result.
\end{proof}
Now, we prove pseudo-stability of \eqref{Rn_morrey} with respect to the gradient-level distance 
$$ A_{\R^N}(u):=\inf_{a,b,x_0}\frac{\|\nabla u-\nabla v_{a,b,x_0}\|_{L^p(\R^N)}}{\| \nabla u\|_{L^p(\R^N)}} \,.$$ 
The proof follows by applying Theorem \ref{Rn_sup}
 and the above reduction lemma.
\begin{thm}\label{Rn-stability} Let $p>N>1.$ Then there exists a positive constant $C=C(p,N)$ such that 
    $$\delta_{M,\mathbb{R}^N}(u)\ge C A_{\R^N}(u^\star)^{p^2}\quad\forall\,u\in W_0^{1,p}(\R^N),$$ for every weakly differentiable function $u$ in $\R^N$ with $|\operatorname{supp}u|<\infty$ and $|\nabla u|\in L^p(\R^N).$
\end{thm}
\begin{proof} Let $\tau(u^\star)$ be as in Theorem \ref{Rn_sup}. If $\tau(u^\star)=0,$ then  $u^\star\in\mathcal{M}_{\mathrm{Morrey}}(\R^N),$ and the theorem follows trivially. Now assume that  $\tau(u^\star)\neq 0$. Notice that $ A_{\R^N}(u^\star)\le 1$. Thus, Theorem \ref{Rn-stability} holds trivially when $\delta_{M,\mathbb{R}^N}(u)\ge 1.$ We therefore assume that $\delta_{M,\mathbb{R}^N}(u)\le 1.$ 
Let $v\in \mathcal{M}_{\mathrm{Morrey}}(\R^N)$ be such that $$\|v\|_{L^\infty(\R^N)}=\|u\|_{L^\infty(\R^N)},\quad \operatorname{supp} u^\star=\operatorname{supp} v, \quad\text{and}\quad\frac{\|u^\star-v\|_{L^\infty(\R^N)}}{\|u^\star\|_{L^\infty(\R^N)}}\le 2\tau(u^\star).$$   Now, 
we apply the Morrey inequality \eqref{Rn_morrey},  Lemma \ref{Thm:Aux:Rn}, and inequality \eqref{Rn_morry-sup} to obtain
 \begin{align*}
     A_{\R^N}(u^\star)^p & \le\left(\frac{\|\nabla u^\star-\nabla v\|_{L^p(\R^N)}}{\| \nabla u^\star\|_{L^p(\R^N)}}\right)^p \le \kappa^p\left(\frac{\|\nabla u^\star-\nabla v\|_{L^p(\R^N)}}{|\operatorname{supp}u^\star|^{\frac{1}{p}-\frac{1}{N}}\| u^\star\|_{L^\infty(\R^N)}}\right)^p\\&\le C[\delta_{M,\mathbb{R}^N}(u)+\frac{\|u^\star-v\|_{L^\infty(\R^N)}}{\|u^\star\|_{L^\infty(\R^N)}}] \ \ \le C[\,\delta_{M,\mathbb{R}^N}(u)+2\,\tau(u^\star)]\\&\le C[\mathbb\,\delta_{M,\mathbb{R}^N}(u)+2C_1\,\delta_{M,\mathbb{R}^N}(u)^\frac{1}{p}] \ \ \le C_2\delta_{M,\mathbb{R}^N}(u)^\frac{1}{p},
 \end{align*}
 for some $C_2=C_2(p,N)>0.$ Taking the power $p$ in the above inequality, we obtain the result.
\end{proof}

\subsection{Gradient estimate for the extremals of the $L^p$-Hardy inequality}
In this subsection, we derive a gradient estimate of the extremals of the $L^p$-Hardy inequality \eqref{Hardy-p}. Our proof is based on standard interior regularity estimates, which plays an important role in our proof of Theorem \ref{Hardy-Thm-grad<2}. We provide a proof of this estimate for completeness.
\begin{lem}\label{regularity}
    Let $\Omega$ be a bounded  $C^{1,\gamma}$ domain in $\mathbb{R}^N$ such that $\Gamma_p(\Omega)>0,$ and let  $\Psi$ be an extremal of the Hardy inequality \eqref{Hardy-p}. Then there exists $C=C(p,N,\Omega)>0$ such that $$|\nabla \Psi(x)|\le  C\frac{\Psi(x)}{d_\Omega(x)}\quad\forall x\in \Omega.$$ 
\end{lem}
\begin{proof} If $\Psi\equiv 0,$ it is trivial. Assume that $\Psi\ne 0.$
 We know that \cite[Section 4,5]{Pinch19} (see also \cite[Theorem 1.1]{Itai}) there exists a constant $C=C(p,N,\Omega)>0$ such that
 \begin{equation}\label{eqv_rl}
     C^{-1}d_\Omega^\alpha(x)\le \Psi(x)\le Cd_\Omega^\alpha(x)\ \ \text{in} \ \Omega
 \end{equation}
for some $\alpha\in \left(1-\frac{1}{p},1\right).$ Thus, it is sufficient to prove $|\nabla\Psi|\le c\, d_\Omega^{\alpha-1}$ in $\Omega,$ for some constant $c>0$. Fix $x_0\in \Omega$ and set $r=d_\Omega(x_0).$ Then $B_r(x_0)\subset \Omega.$ Define a scaling function $$\phi_r(y):=\frac{\Psi(x_0+ry)}{r^\alpha}\,,\quad y\in B_{\frac{1}{2}}(0).$$ Since $d_\Omega$ is a Lipschitz function and $r=d_\Omega(x_0)$, 
\begin{equation}\label{Eq2:ap2}
    \frac{r}{2} \le d_\Omega(x_0+ry)\le \frac{3r}{2} \quad\forall\,y\in B_{\frac{1}{2}}(0).
\end{equation}
Now, using \eqref{eqv_rl} and \eqref{Eq2:ap2} we infer that 
\begin{equation}\label{eq3;ap2}
    \|\phi_r\|_{L^\infty(B_{\frac{1}{2}}(0))}\le C\frac{d_\Omega(x_0+ry)^\alpha}{r^\alpha}\le 4^\alpha C\le 4C.
\end{equation}
Since $\Psi$ satisfies \eqref{EVP}, 
one can verify that $\phi_r$ satisfies the following PDE:
 $$-\Delta_p \phi_r = \mathbb{H}_p(\Omega)^p\left(\frac{r}{d_\Omega(x_0+r \ \cdot)}\right)^p\phi_r^{p-1} =:f_r\quad\text{in}\,B_{\frac{1}{2}}(0).$$ Moreover, combine \eqref{Eq2:ap2} and \eqref{eq3;ap2} to get $\|f_r\|_{L^\infty(B_{\frac{1}{2}}(0))}\le C_1,$ for some $C_1=C_1(p,N,\Omega)>0$. 
   By the standard interior $C^{1,\beta},\,\beta\in (0,1)$-regularity theory for the $p$-Laplace equation \cite{Di}, $$\|\nabla \phi_r\|_{L^\infty(B_{\frac{1}{16}}(0))}\le C_2$$ for some  $C_2=C_2(p,N,\Omega)>0.$ 
   Evaluating at $y=0$, $$|\nabla \Psi(x_0)|=r^{\alpha-1}|\nabla \phi_r(0)|\le C_2 r^{\alpha-1}=C_2d_\Omega(x_0)^{\alpha-1}.$$ Since $x_0$ is arbitrary, we conclude the result.
\end{proof}

\subsection{Other auxiliary results}

\begin{thm} \label{Thm:A}
    Let $(X,\|\cdot\|)$ be a normed linear space and $Y$ be a one-dimensional subspace of $X$. Then $$D(x) \leq D_{\mathrm{c}}(x) \leq 2 D(x) \,, \ \ \forall x \in X \,,$$ 
    where $D(x):=\inf\{\|x-y\| : y \in Y\}$ is the usual distance from $x$ to $Y$ and $D_c(x):=\inf\{\|x-y\| : y \in Y \  \mbox{with} \ \|y\|=\|x\|\}$ is the constrained version on $D$.
\end{thm}
\begin{proof}
It is obvious that  $D(x) \leq D_{\mathrm{c}}(x)$ for all $x \in X$. We will show that  $D_{\mathrm{c}}(x) \leq 2 D(x)$ for all $x \in X$. Let $Y={\mathrm{Span}}(e)$, where $e$ is a unit vector in $X$. Clearly, 
\begin{align*}
 D(x) &=\inf \{\|x-ce\|: c \in \mathbb{R}\} \,, \\
D_c(x) &=\min \{\|x-\|x\|e\|,  \|x+\|x\|e\|\} \,.
\end{align*}
Since $Y$ is one-dimensional, $D(x)$ is attained \cite[Chap. 1, Corollary 2.2]{Singer}, i.e., there exists $c_x \in \mathbb{R}$ such that
$$D(x)=\|x-c_x e\| \,.$$
Now, the triangular inequity gives
    \begin{equation*}
    \|x\|\le \|x-c_x e\| + \|c_x e\|=D(x)+|c_x| \,,
    \end{equation*}
 and
$$|c_x|=\|c_x e\| \le \|x-c_x e\| + \|x\|=D(x)+\|x\| \,.$$ 
Together, we have 
    \begin{equation}\label{11:ste2}
       ||c_x|-\|x\||\le D(x). 
    \end{equation}
If $c_x \geq 0,$ then  $$D_c(x) \le \|x - \|x\| e\| \le \| x - c_x e \| + \|(c_x -\|x\|)e\| = \| x - c_x e \| + |\,|c_x| -\|x\||  \le 2D(x).$$ 
Conversely, if
$c_x < 0,$ then  $$D_c(x) \le \|x + \|x\| e\| \le \| x - c_x e \| + \|(c_x +\|x\|)e\| = \| x - c_x e \| + ||c_x| -\|x\||  \le 2D(x).$$
Combining both cases, we conclude the theorem.
\end{proof}

\begin{center}
	{\bf Acknowledgement}
\end{center}
 S. Roy, acknowledges the support of the National Board for Higher Mathematics (NBHM) through the NBHM Postdoctoral Fellowship
 (No.0202/3/2026/R\&D-II/2870). The authors thank Prof. Anup Biswas for the discussion on the gradient estimate for the extremals of the Hardy inequality. The authors also thank Prof. Yehuda Pinchover for his valuable remarks on the article.

\bibliographystyle{abbrvurl}
\bibliography{Reference}

\end{document}